\documentclass[10pt]{scrartcl}

\usepackage[a4paper, left=1.5cm, right=1.5cm, top=2cm]{geometry}

\usepackage[T1]{fontenc}
\usepackage{inputenc}
\usepackage{algorithm}
\usepackage{amsthm}
\usepackage{amssymb}
\usepackage{algpseudocode} 
\usepackage{wrapfig}
\usepackage{graphicx}
\usepackage{float}
\usepackage{mathtools}
\mathtoolsset{showonlyrefs} 
\usepackage{amsmath} 
\usepackage{caption}
\numberwithin{equation}{section} 
\usepackage{cases}
\usepackage{tabularx}
\usepackage{xcolor}
\usepackage{url}
\usepackage{enumitem}
\usepackage{authblk}
\usepackage{lmodern}

\newtheorem{theorem}{Theorem}
\newtheorem{lemma}{Lemma}

\theoremstyle{remark}
\newtheorem{remark}{Remark}
\theoremstyle{definition}

\newtheorem{definition}{Definition}
\newtheorem{problem}{Problem}
\newtheorem{assumption}{Assumption}

\def\orcid[#1]#2#3{}%
\newcommand{\orgdiv}[1]{#1}%
\newcommand{\orgname}[1]{#1}%
\newcommand{\orgaddress}[1]{#1}%
\newcommand{\state}[1]{#1}%
\newcommand{\country}[1]{#1}%

\def\R{{\mathbb{R}}}
\def\N{{\mathbb{N}}}
\newcommand{\ti}{\tilde}
\renewcommand{\S}{\mathcal{S}}

\newcommand{\h}{\hat}

\DeclareMathOperator*{\esssup}{ess\,sup}

\title{A Gauss-Newton Method for ODE Optimal Tracking Control}

\author[1]{Vicky Holfeld}

\author[1]{Michael Burger}

\author[2]{Claudia Schillings}

\affil[1]{\normalsize \orgdiv{Division Mathematics for Vehicle Engineering}, \orgname{Fraunhofer Institute for Industrial Mathematics (ITWM)}, \orgaddress{\state{Kaiserslautern}, \country{Germany}}}

\affil[2]{\normalsize \orgdiv{Research Group Numerical Analysis of Partial Differential Equations}, \orgname{Free University of Berlin}, \orgaddress{\state{Berlin},\\ \country{Germany}}}

\date{\large \today}

\begin{document}

\maketitle

\begin{center}
\end{center}
\vspace{0.5cm}

\begin{center}
	\textbf{Abstract}
\end{center}
\begin{abstract}
  This paper introduces and analyses a continuous optimization approach to solve optimal control problems involving ordinary differential equations (ODEs) and tracking type objectives. Our aim is to determine control or input functions, and potentially uncertain model parameters, for a dynamical system described by an ODE. We establish the mathematical framework and define the optimal control problem with a tracking functional, incorporating regularization terms and box-constraints for model parameters and input functions. Treating the problem as an infinite-dimensional optimization problem, we employ a Gauss-Newton method within a suitable function space framework. This leads to an iterative process where, at each step, we solve a linearization of the problem by considering a linear surrogate model around the current solution estimate. The resulting linear auxiliary problem resembles a linear-quadratic ODE optimal tracking control problem, which we tackle using either a gradient descent method in function spaces or a Riccati-based approach. Finally, we present and analyze the efficacy of our method through numerical experiments.
\end{abstract}

\section{Introduction }\label{sec1} 
Optimal control generally deals with the problem of finding a control input function for a dynamical system subject to constraints and with the objective of optimizing a certain performance criterion. It has a wide range of applications and can be found in nearly any discipline like, for instance, in natural sciences, engineering sciences, economy and biology, see \cite{Gerdts2012}, \cite{Betts2020}, \cite{Bryson}.
A specific subclass hereof is optimal tracking control. Problems therein, also referred to as tracking problems, cf. \cite{Sontag1998}, \cite{Locatelli2001}, aim to compute controls that reproduce a given target or reference output for the dynamical system. In this contribution, we particularly focus on such problems, especially, on continuous-time optimal tracking control problems subject to ordinary differential equations (ODEs), where an ODE as a constraint describes the behavior of the dynamical system over time.

Tracking problems find applications across various domains. In the engineering sciences, for instance, tracking reference dynamics is used for the optimization of aircraft (see, e.g., \cite{AircraftTracking}) and spacecraft trajectories  (see, e.g., \cite{Betts2020}, \cite{SpacecraftTracking}) but also for controlling wafer-stages for the production of computer-chips \cite{Dijkstra2004}.
Moreover, tracking problems often arise in optimizing the path of mobile robots \cite{Faulwasser2012} and vehicles \cite{AutonomousVehiclesTracking}. Optimal tracking control is also common in vehicle engineering in the context of testing and vehicle simulation \cite{Burger2014}, \cite{ThesisBurger2011}, 
which will be in our focus in this contribution.

Tracking problems are often not solvable analytically and, thus, numerical solutions to such problems are usually sought in practice. So far, there exist a large number of numerical methods for general time-continuous ODE optimal control problems, that are also applicable to the subclass of tracking problems. For an overview, we refer to \cite{Betts2020}, \cite{Bryson} and, in particular, \cite{Gerdts2012}, where the more general setting of optimal control for differential algebraic equations (DAEs) is considered. These methods can be classified as \textit{first-discretize-then-optimize} and \textit{first-optimize-then-discretize} methods. For numerical approaches in the first class, also known as direct methods, the time-continuous ODE optimal control problem is a priori transformed to a time-discrete problem via a discretization approach. The resulting, typically nonlinear, optimization problem can then be solved by techniques and methods from nonlinear programming, cf. \cite{Bertsekas2008}, \cite{NocedalWright}. For methods in the second class, also referred to as function space approaches, the time-continuous optimal control problem is treated as an infinite-dimensional optimization problem and is solved in an appropriate Banach or Hilbert space setting. Besides of the better-known indirect solution methods, that are based on solving first-order necessary optimality conditions, the more rarely appearing direct function space methods are derived by extending well-known algorithms for finite-dimensional optimization problems to the Banach or Hilbert space setting and applying them directly to the optimal control problem. An overview of these function space methods for ODE optimal control problems can be found in \cite{Polak1972}. For example, a function space gradient descent method, see \cite{Gerdts2012}, relies on the computation of a descent direction via the corresponding adjoint equations. For control-constrained problems, a function space Lagrange-Newton method can, for instance, be found in \cite{Alt1993} and a sequential quadratic programming (SQP) method for continuous-time ODE optimal control problems is presented in \cite{Machielsen1988}.

In this contribution we aim to derive and apply a function space method that is specifically applicable to ODE optimal tracking control problems and has certain advantages in practical implementation in terms of efficiency and applicability to real-world problems. We consider tracking problems, where we aim to compute simultaneously control inputs as well as unknown model parameters for the dynamical system and for which box-constraints for controls and parameters are additionally allowed to be included. Treating these as infinite-dimensional optimization problems and exploiting the special structure, namely nonlinear least-squares problem, we directly apply a Gauss-Newton method, generally formulated in Banach spaces (see, e.g., \cite{Fu2020}, \cite{Kaltenbacher2009}, \cite{Kaltenbacher2010} for the iteratively regularized Gauss-Newton method for ill-posed nonlinear operator equations in Banach spaces) or Hilbert spaces (see, e.g., \cite{Kaltenbacher2008}, \cite{Amat2014}), to our continuous-time problem and derive a projected Gauss-Newton function space method for the considered tracking problems. 

The remaining part of the article is organized as follows. In Section 2, we define dynamical systems by parameter-dependent and controlled ODEs and present the continuous-time ODE optimal control problem with tracking objective and simple box-constraints. Then, in Section 3, we derive and formulate the projected Gauss-Newton approach in function spaces as a new solution method for the aforementioned problem and formulate a global convergence result in a specific setting. In Section 4, we show a numerical example where the presented method is applied and finally compare the method to the gradient descent method in function spaces to show the efficiency of our proposed method. Section 5 concludes the article with a summary and short discussion. \\

\noindent
\textit{Notation.} Throughout this paper, we denote with $\Vert \cdot \Vert_{\R^{n}}$ the Euclidean norm in $\R^{n}$ with corresponding dimension $n \in \N.$ 
For an interval $I := [t_0, t_f] \subseteq \R$, $t_0 < t_f$, we denote by $L^\infty(I;\R^{n})$ the Banach space of all Lebesgue measurable functions $x: I \rightarrow \R^{n}$ that are essentially bounded, i.e., $\| x\|_{\infty} := \| x\|_{L^\infty(I;\R^n)} := \esssup\limits_{t \in I} \| x(t) \|_{\R^{n}} := \underset{N \subset I, ~ \mu(N) = 0}{\inf} ~ \underset{t \in I \setminus N}{\sup} \| x(t) \|_{\R^{n}}$ is finite. We further denote by $W^{1,\infty}(I,\R^{n})$ the Banach space of all absolutely continuous functions $x: I \rightarrow \R^{n}$ with $\Vert x \Vert_{1,\infty} < \infty$, where $\Vert x \Vert_{1,\infty} := \max \{\Vert x \Vert_{\infty}, \Vert \dot{x} \Vert_{\infty}\}$. The space $L^{2}(I;\R^{n})$ denotes the Hilbert space of all Lebesgue measurable functions $x: I \rightarrow \R^{n}$ satisfying $ \| x\|_{L^{2}(I;\R^{n})} < \infty$, with norm $ \| x\|_{L^{2}(I;\R^{n})} := \sqrt{\left\langle x, x\right\rangle_{L^{2}(I;\R^{n})}} $ induced by the inner product $\left\langle x, \ti x\right\rangle_{L^{2}(I;\R^{n})} := \int_{t_0}^{t_f} x(t)^\top \ti x(t) ~ dt$. We use a.e. and a.a. as abbreviations for almost everywhere and almost all. 

\section{The Constrained and Parameter-Dependent ODE Optimal Control Problem with Tracking Functional and Preliminaries}\label{sec2}
In this section, we introduce the mathematical framework for the  formulation of the tracking problem and for the derivation of a novel numerical solution method for these types of problems.

\subsection{The Dynamical System}\label{sec21} 
\noindent
We consider a dynamical system whose state, denoted by $x$, is mathematically described by a parameter-dependent and controlled ODE initial value problem (ODE-IVP)
\begin{equation}\label{DynSys}
	\begin{cases}
		\dot{x}(t) &= f(t,x(t),u(t),p) ~~ \mbox{a.e. in }I = [t_0, t_f]\\
		x(t_0)&= x_0.          
	\end{cases}
\end{equation}
with given initial state $x_0 \in \R^{n_x}$, $n_x \in \N$, and a parameter vector $p \in \R^{n_p}$, $n_p \in \N$ as well as an external control function $u \in L^{\infty}(I;\R^{n_u})$, $n_u \in \N$, as system inputs. In the sequel, we want to ensure that the ODE-IVP has a unique solution. Therefore, we make the following assumption. 
\noindent
\begin{assumption}
\label{Ass1} The right-hand side $f: I \times \R^{n_x} \times \R^{n_u} \times \R^{n_p} \longrightarrow \R^{n_x}$  in \eqref{DynSys} is assumed to be twice continuously differentiable and globally Lipschitz continuous. The latter means that there is a constant $L \geq 0$ such that
\[\| f(t,x,u,p)-f(t,\bar{x},\bar{u},\bar{p}) \|_{\R^{n_x}} \leq L \left( \| x-\bar{x}\|_{\R^{n_x}} + \| u-\bar{u}\|_{\R^{n_u}} + \| p-\bar{p} \|_{\R^{n_p}} \right) \]
for all $t \in I $ and all $x,\bar{x} \in \R^{n_x}$, $u,\bar{u} \in \R^{n_u}$, $p, \bar{p} \in \R^{n_p}$.  \\
\end{assumption}

\noindent
As a consequence of Assumption\ref{Ass1}, the ODE-IVP \eqref{DynSys} possesses a unique global solution $x \in W^{1,\infty}(I;\R^{n_x})$ for any pair $(u,p) \in L^{\infty}(I;\R^{n_u}) \times \R^{n_p}$, see, e.g., \cite{Sontag1998} for ODE solutions in the sense of Carath\'eodory. Having a unique solution, we are able to define a solution operator for \eqref{DynSys} that maps the system inputs to the system states. \\

\begin{definition}
Let a control function $u \in L^{\infty}(I;\R^{n_u})$ and a parameter vector $p \in \R^{n_p}$ be given. We define the solution operator $\ti{\S}$, also referred to as the \textit{input-to-state map}, that assigns each pair $(u,p)$ the unique global solution $x \in W^{1,\infty}(I; \R^{n_x})$ of the aforementioned ODE-IVP \eqref{DynSys}, i.e.
\begin{equation}
	\ti{\S}: L^{\infty}(I;\R^{n_u}) \times \R^{n_p} \longrightarrow W^{1,\infty}(I; \R^{n_x}),  \quad \ti{\S}(u,p) = x.
\end{equation}\\
\end{definition} 

\noindent
Furthermore, we extend the ODE-IVP by an \textit{output-equation}, that is for $t \in I$, 
\begin{align}\label{OutputEqu}
	y(t) = h(t,x(t),u(t),p),
\end{align}
with outputs $y$ and nonlinear function $h$ depending on the state $x$ and, possibly, also of the control function $u$ and the parameter vector $p$. 
We make the following assumptions about the differentiablity. 

\begin{assumption} \label{Ass2} The output-function $h: I \times \R^{n_x} \times \R^{n_u} \times \R^{n_p} \longrightarrow \R^{n_y}$ in \eqref{OutputEqu}
is assumed to be twice continuously differentiable w.r.t. all arguments and, in addition, $n_y \leq n_x$ is assumed. \\
\end{assumption}

\noindent
We define the operator $\S$ describing the dynamical system given by the state-space representation \eqref{DynSys}-\eqref{OutputEqu}, that maps the system inputs $(u,p)$ to the system outputs $y$. 

\noindent
\begin{definition} Similarly to the input-to-state map, we define the \textit{input-to-output map} $\S$ for a given control function $u \in L^{\infty}(I;\R^{n_u})$ and parameter vector $p \in \R^{n_p}$ as 
\begin{equation}
	\begin{split}
		& \S: L^{\infty}(I;\R^{n_u}) \times \R^{n_p} \longrightarrow L^{\infty}(I;\R^{n_y}), \quad \left[\S(u,p)\right](t) = h(t,x(t),u(t),p),
	\end{split}
\end{equation}
with $x = \ti{\S}(u,p)$. \\
\end{definition}

\noindent
Note that the Assumptions \ref{Ass1}, \ref{Ass2} are required for the differentiablity of the operator $\S$ that will be presented in Subsection \ref{sec23}. 

With this state-space representation for the dynamical system, represented by the operator $\S$, we are now able to mathematically formulate the tracking problem. 

\subsection{Problem Formulation - The ODE Optimal Control Problem with Tracking Functional}\label{sec22}
In the following, we are interested in the inverse problem, that is, to find system inputs $(u,p) \in L^{\infty}(I;\R^{n_u}) \times \R^{n_p}$ for the dynamical system $\S : L^{\infty}(I;\R^{n_u}) \times \R^{n_p} \rightarrow L^{\infty}(I;\R^{n_y})$ that reproduce a given reference output $y_\text{ref} \in L^{\infty}(I;\R^{n_y})$ in the best possible way. In addition, we introduce some box-constraints on the system inputs $u$ and $p$ to restrict their values in a certain range. Therefore, we define the closed, convex set 
\begin{align}\label{ConvexSetConstraints}
	\mathcal{U} := \left\lbrace (u,p) \in L^{\infty}(I;\R^{n_u}) \times \R^{n_p} ~\big|~ u^\text{low}_i\leq u_i(t)\leq u^\text{up}_i \mbox{ a.e. in }I,~ p^\text{low}_j\leq p\leq p^\text{up}_j, ~i=1,..,n_u,~	j=1,..,n_p\right\rbrace,
\end{align}
$\mathcal{U} \subseteq L^{\infty}(I;\R^{n_u}) \times \R^{n_p}$, with fixed bounds $u^\text{low}_i, u^\text{up}_i, p^\text{low}_j, p^\text{up}_j \in \R$.
Hence, we formulate the inversion task by the following constrained Tikhonov-regularized optimization problem. 

\begin{problem}
\label{ODEOptControlProblem}(\textit{Nonlinear Least-Squares Problem}) 
Find system inputs $(u,p) \in \mathcal{U} \subseteq  L^\infty(I;\R^{n_u}) \times \R^{n_p}$ such that the error functional $J: \mathcal{U} \rightarrow \R$ given by
\begin{align}\label{ODEOptControl}
	J(u,p) &= \frac{1}{2} \left\Vert Q^{\frac{1}{2}} \left[\S(u,p) - y_\text{ref}\right] \right\Vert_{L^2(I;\R^{n_y})}^2 + \frac{\alpha_u}{2} \left\Vert u \right\Vert_{L^2(I;\R^{n_u})}^2 + \frac{\alpha_p}{2} \left\Vert p \right\Vert_{\R^{n_p}}^2 + \frac{1}{2} \left\Vert T^{\frac{1}{2}} \left[\S(u,p)(t_f) - y_\text{ref}(t_f)\right] \right\Vert_{R^{n_y}}^2 
\end{align}
is minimized for regularization parameters $\alpha_u, \alpha_p > 0$, symmetric positive semi-definite weight matrices $Q, T \in \R^{n_y \times n_y}$ and a given reference $y_\text{ref} \in L^\infty(I;\R^{n_y})$. \\
\end{problem}

By means of the input-to-output operator $\S$ defined in Subsection \ref{sec21} and the use of the $L^2$-norm, Problem \ref{ODEOptControlProblem} has the form of a nonlinear ODE optimal control problem with tracking objective function and box-constraints for the control and the parameter vector. This is stated in the following equivalent problem formulation. 

\begin{problem}
\label{NonlinearODEOCP}(\textit{ODE Optimal Control Problem with Tracking Functional}) Find a control function $u \in L^{\infty}(I;\R^{n_u})$ and a parameter vector $p \in \R^{n_p}$ that minimize 
\begin{equation}\label{ODE Optimal Control OptProblem}
	\begin{split}
		J(u,p) &=  \int_{t_0}^{t_f} \frac{1}{2} \left\Vert Q^{\frac{1}{2}} \left[h(t,x(t),u(t),p)-y_\text{ref}(t)\right] \right\Vert_{\R^{n_y}}^2 + \frac{\alpha_u}{2} \left\Vert u(t) \right\Vert_{\R^{n_u}}^2 ~dt + \frac{\alpha_p}{2} \left\Vert p \right\Vert_{\R^{n_p}}^2 \\
		&~~~~ + \frac{1}{2} \left\Vert T^{\frac{1}{2}} \left[h(t_f,x(t_f),u(t_f),p) - y_\text{ref}(t_f)\right] \right\Vert_{R^{n_y}}^2  \\
	\end{split}	
\end{equation}
such that 
\begin{equation}
	\begin{split}
		&  \begin{cases}
			\dot{x}(t) &= f(t,x(t),u(t),p), ~~ t \in I, \\
			x(t_0)&= x_0,
		\end{cases}
	\end{split}	
\end{equation}
and such that the following constraints are satisfied
\begin{equation}
	\begin{split}
		& u^\text{low}_i \leq u_i(t) \leq u^\text{up}_i ~~ \mbox{for a.a.} ~t \in  I, ~~ i = 1,...,n_u, \\
		& p^\text{low}_j \leq p_j \leq p^\text{up}_j, ~~ j = 1,...,n_p. \\
	\end{split}
\end{equation}\\
\end{problem} 

\begin{remark}
	Here, we focus on solving the inverse problem 
	by considering the (constrained) Tikhonov-regularized optimization problem \ref{ODEOptControlProblem} for $(u,p)$ for given regularization parameters $\alpha_u, \alpha_p$. The (optimal) choice of the regularization parameters $\alpha_u, \alpha_p$ can be motivated via the Bayesian approach, in particular the maximum a-posteriori estimate . We refer to \cite{Kaipio2005}, \cite{Dashti2017}, \cite{Stuart2010}  for more details. \\
\end{remark}

\subsection{The Fr\'echet-Derivative}\label{sec23}
We aim to derive a novel numerical method that solves the ODE optimal control problem \ref{NonlinearODEOCP} iteratively by subproblems for which we require the derivative of the error function $J$ and the input-to-output map $\S$, respectively. Since $S$ and $J$ are operators in Banach spaces, the concept of Fr\'echet-derivatives is used in the following, cf. \cite{Gerdts2012}. 

For the Fr\'echet-derivative of $J$ and $\S$ we first establish the Fr\'echet-differentiablity of the operator $\ti{\S}$. 

\begin{theorem}\label{Theorem_Fréchet}
	The input-to-state map $\ti{\S}$ is continuously Fr\'echet-differentiable. In particular, the Fr\'echet-derivative of $\tilde{\S}$ at a point $(\hat{u},\hat{p})$ with $\hat{x} = \ti{\S}(\hat{u},\hat{p})$  is given by
	\begin{equation}
		\partial_{(u,p)}\ti{\S} \big|_{(\hat{u},\hat{p})}: L^{\infty}(I;\R^{n_u}) \times \R^{n_p} \longrightarrow W^{1,\infty}(I; \R^{n_x})
	\end{equation}
	\begin{equation}
		\partial_{(u,p)}\ti{\S} \big|_{(\hat{u},\hat{p})}(\delta u,\delta p) = \delta x
	\end{equation}
	with $\delta x$ being the unique global solution of the linear time-varying ODE-IVP
	\begin{equation}\label{eq:FrechDiff}
		\begin{split}
			\dot{\delta x}(t) &= f_x[t]\delta x(t) + f_u[t]\delta u(t) + f_p[t]\delta p, \\
			\delta x (t_0) &= 0.
		\end{split}
	\end{equation}
	For readability, we make use of the short-hand notations $f_x[t] := f_x(t,\h{x}(t),\h{u}(t),\h{p})$, $f_u[t]:= f_u(t,\h{x}(t),\h{u}(t),\h{p})$ and $f_p[t]:= f_p(t,\h{x}(t),\h{u}(t),\h{p})$.
\end{theorem}

\noindent
The proof of Theorem \ref{Theorem_Fréchet} can be found in \cite{ThesisBurger2011}. There, it is shown for the more general case of a differential algebraic equation (DAE) as model equations for the dynamical system $\S$. 
Next, We obtain the differentiablity of the input-to-output-map $\S$ by the following. 
\begin{lemma}\label{Lemma1}
	The input-to-output map $\S$ is continuously Fr\'echet-differentiable and the Fr\'echet-derivative of $\S$ at a point $(\hat{u},\hat{p})$ is given by 
	\begin{equation}
		\partial_{(u,p)}\S \big|_{(\hat{u},\hat{p})}: L^{\infty}(I;\R^{n_u}) \times \R^{n_p} \longrightarrow L^{\infty}(I; \R^{n_y})
	\end{equation}
	\begin{equation}
		\left[ \partial_{(u,p)}\S \big|_{(\hat{u},\hat{p})}(\delta u,\delta p) \right] (t) =  h_x[t]\delta x(t) + h_u[t]\delta u(t) + h_p[t]\delta p,
	\end{equation}
	where $h_x[t] := h_x(t,\h{x}(t),\h{u}(t),\h{p})$ for $\hat{x} = \ti{\S}(\hat{u},\hat{p})$. The other Jacobians $h_u[t], h_p[t]$ are defined correspondingly and $\delta x = \partial_{(u,p)}\ti{\S} \big|_{(\hat{u},\hat{p})}(\delta u,\delta p)$. For simplification, we write $S^\prime(\hat{u},\hat{p})$ for $\partial_{(u,p)} S \big|_{(\hat{u},\hat{p})}$ in the following.
\end{lemma}
\noindent
The proof of the preceding lemma directly follows by the chain-rule for Fr\'echet-derivatives, cf. \cite{IoffeTihomirov}, and the Fr\'echet-differentiability of the input-to-state map $\ti \S$ stated in Theorem \ref{Theorem_Fréchet}.

For the Fr\'echet-derivative of the cost functional $J: L^{\infty}(I;\R^{n_u}) \times \R^{n_p} \longrightarrow \R$, we exploit that $J$, in \eqref{ODE Optimal Control OptProblem}, can in general be written as
\begin{align*}
	J(u,p) = \phi(t_f,x(t_f),u(t_f),p) + \int_{t_0}^{t_f} \varphi(t,x(t),u(t),p) ~ dt
\end{align*} 
with functions $\phi, \varphi : \R \times \R^{n_x} \times \R^{n_u} \times \R^{n_p} \longrightarrow \R$ being given as
\begin{align}
	\phi(t_f,x(t_f),u(t_f),p) &= \frac{1}{2} \left\Vert T^{\frac{1}{2}} \left[h(t_f,x(t_f),u(t_f),p) - y_\text{ref}(t_f)\right] \right\Vert_{R^{n_y}}^2, \\
	\varphi(t,x(t),u(t),p) &= \frac{1}{2} \left\Vert Q^{\frac{1}{2}} \left[h(t,x(t),u(t),p)-y_\text{ref}(t)\right] \right\Vert_{\R^{n_y}}^2 + \frac{\alpha_u}{2} u(t)^\top u(t) + \frac{\alpha_p}{2(t_f - t_0)}p^\top p.
\end{align}
We then get the differentiability of the cost function $J$ and its representation of the Fr\'echet-derivative as follows.
\begin{lemma}\label{Lemma2}
	The error function $J$ is continuously Fr\'echet-differentiable. The Fr\'echet-derivative of $J$ at a point $(\hat{u},\hat{p})$ is given by 
	\begin{align}
		\partial_{(u,p)} J \big|_{(\hat{u},\hat{p})}: L^{\infty}(I;\R^{n_u}) \times \R^{n_p} \longrightarrow \R
	\end{align}
	with
	\begin{align}
		\partial_{(u,p)} J\big|_{(\hat{u},\hat{p})}(\delta u, \delta p) = \phi_x[t_f] \delta x(t_f) + \phi_u[t_f] \delta u(t_f) + \phi_p[t_f] \delta p + \int_{t_0}^{t_f} \varphi_x[t] \delta x(t) + \varphi_u[t] \delta u(t) + \varphi_p[t] \delta p ~ dt,
	\end{align}
	where $\delta x = \partial_{(u,p)} \ti \S \big|_{(\hat{u},\hat{p})}(\delta u,\delta p)$, $\varphi_x[t] = \varphi(t, \hat{x}(t),\hat{u}(t),\hat{p})$ for $\hat{x} = \ti \S (\hat{u},\hat{p})$ and similar for the other Jacobians of $\phi$ and $\varphi$. For a simplified notation, we write $J^\prime(\hat{u},\hat{p})$ for $\partial_{(u,p)} J \big|_{(\hat{u},\hat{p})}$ in the following.
\end{lemma}
\noindent
For the proof one shows that, under Assumption \ref{Ass2}, $\phi$ and $\varphi$ are twice continuously differentiable and, again, applies the chain-rule and the Fr\'echet-differentiability of the input-to-state map $\ti \S$ to get the desired result. 

\section{The Gauss-Newton Method}\label{sec3}
In this section, the goal is to derive a 
solution method for the ODE optimal control problem \ref{ODEOptControlProblem} and to provide a global convergence result.

\subsection{Derivation of the Gauss-Newton procedure}
Recall the aim is to solve Problem \ref{ODEOptControlProblem} that has the form of a nonlinear least-squares problem. Assuming that $\mathcal{U} = L^{\infty}(I;\R^{n_u}) \times \R^{n_p}$, we solve the (unconstrained) problem numerically by applying the Gauss-Newton method in $L^\infty(I;\R^{n_u}) \times \R^{n_p}$. This results in an iterative procedure of the form 
\begin{align}\label{IterProc}
	\left(u^{(k+1)},p^{(k+1)}\right) = \left(u^{(k)}, p^{(k)}\right) + \gamma_k \cdot \left(\delta u^{(k)}, \delta p^{(k)}\right), \quad k = 0,1,2,...
\end{align}
with stepsize $\gamma_k > 0$.
The updates $(\delta u^{(k)}, \delta p^{(k)}) \in L^\infty(I;\R^{n_u}) \times \R^{n_p}$ are determined as solutions of subproblems that are achieved by a linearization of the objective function around the current iterate $(u^{(k)},p^{(k)})$, given by a first-order Taylor approximation, that is minimized in the least-squares sense. In our setting, the subproblems are of the following form.
\begin{problem}
\label{LinearLeastProb} (\textit{Linear Least-Squares Problem}) Find $\delta u^{(k)} \in L^{\infty}(I;\R^{n_u})$ and $\delta p^{(k)} \in \R^{n_p}$ that minimize
\begin{equation}\label{AuxiliaryProblem_J}
	\begin{split}
		\h J_\alpha (\delta u^{(k)}, \delta p^{(k)}) &:= \frac{1}{2} \left\Vert Q^{\frac{1}{2}}\left[\S(u^{(k)},p^{(k)}) + \S^\prime (u^{(k)},p^{(k)})(\delta u^{(k)}, \delta p^{(k)}) - y_\text{ref}\right] \right\Vert^2_{L^2(I;\R^{n_y})} + \frac{\alpha_u}{2} \left\Vert u^{(k)} + \delta u^{(k)} \right\Vert^2_{L^2(I;\R^{n_u})} \\ 
		& ~~~~~~~ + \frac{\alpha_p}{2} \left\Vert p^{(k)} + \delta p^{(k)} \right\Vert^2_{\R^{n_p}}  +  \frac{1}{2} \left\Vert T^{\frac{1}{2}}\left[\S(u^{(k)},p^{(k)})(t_f) + \S^\prime (u^{(k)},p^{(k)})(\delta u^{(k)}, \delta p^{(k)})(t_f) - y_\text{ref}(t_f)\right] \right\Vert^2_{\R^{n_y}}.  \\	
	\end{split}
\end{equation} \\
\end{problem}

\noindent
Using the definition of $\S$ and its Fr\'echet-derivative $\S ^\prime(u,p)$, see Lemma \ref{Lemma1}, Problem \ref{LinearLeastProb} has the form of a linear-quadratic ODE optimal control problem in $(\delta u^{(k)}, \delta p^{(k)})$ as stated in Problem \ref{LQODEOCP}. Note, that, compared to Problem $\ref{LinearLeastProb}$, here we skip the terms that are constant with respect to $(\delta u^{(k)}, \delta p^{(k)})$. 

\begin{problem}
\label{LQODEOCP}(\textit{Linear-Quadratic ODE Optimal Control Problem}) Find a control function $\delta u^{(k)} \in L^{\infty}(I;\R^{n_u})$ and a parameter vector $\delta p^{(k)} \in \R^{n_p}$ that minimize
{\small
\begin{equation}\label{AuxiliaryProblem}
	\begin{split}
		\hat{J}_{\alpha}\left(\delta u^{(k)}, \delta p^{(k)}\right)	&=  \int_{t_0}^{t_f} \frac{1}{2} \left\Vert Q^{\frac{1}{2}}\left[ C(t)\delta x(t) +  D_u(t)\delta u^{(k)}(t) + D_p(t)\delta p^{(k)} - r(t) \right] \right\Vert_{\R^{n_y}}^2  + \frac{\alpha_u}{2} \left[ 2 u^{(k)}(t)^\top \delta u^{(k)}(t) + \left\Vert \delta u^{(k)}(t) \right\Vert_{\R^{n_u}}^2\right] dt \\
		& ~~~ + \frac{\alpha_p}{2} \left[ 2 {p^{(k)}}^\top \delta p^{(k)} + \left\Vert \delta p^{(k)} \right\Vert_{\R^{n_p}}^2 \right] +  \frac{1}{2} \left\Vert T^{\frac{1}{2}}\left[C(t_f)\delta x(t_f) +  D_u(t_f)\delta u^{(k)}(t_f) + D_p(t_f)\delta p^{(k)} - r(t_f)\right] \right\Vert^2_{\R^{n_y}}
	\end{split}  
\end{equation}}
such that 
\begin{equation}\label{AuxiliaryProb_ODEConstraints}
	\begin{cases}
		\dot{\delta x}(t) &= A(t)\delta x(t) + B_u(t)\delta u^{(k)}(t) + B_p(t) \delta p^{(k)}, ~~ t \in I, \\
		\delta x(t_0)&= 0,         
	\end{cases}	
\end{equation}
with matrix functions 
\begin{align*}
	A(t) &= f_x(t,x^{(k)}(t),u^{(k)}(t),p^{(k)}) \in \R^{n_x \times n_x}, \quad B_u(t) = f_u(t,x^{(k)}(t),u^{(k)}(t),p^{(k)}) \in \R^{n_x \times n_u}, \\
	B_p(t) &= f_p(t,x^{(k)}(t),u^{(k)}(t),p^{(k)}) \in \R^{n_x \times n_p}, \quad C(t) = h_x(t,x^{(k)}(t),u^{(k)}(t),p^{(k)}) \in \R^{n_y \times n_x}, \\
	D_u(t) &= h_u(t,x^{(k)}(t),u^{(k)}(t),p^{(k)}) \in \R^{n_y \times n_u}, \quad D_p(t) = h_p(t,x^{(k)}(t),u^{(k)}(t),p^{(k)}) \in \R^{n_y \times n_p},
\end{align*} 
$x^{(k)} = \ti \S \left(u^{(k)},p^{(k)}\right)$ and  $r(t) := y_\text{ref}(t) - \S\left(u^{(k)},p^{(k)}\right)(t)$. \\
\end{problem} 

\begin{remark} (\textit{Existence and Uniqueness of a Solution for the Auxiliary Problem}) Since $\hat{J}_{\alpha}$ is Fr\'echet-differentiable, it is continuous and due to Weierstra\ss~ theorem (see, for instance, \cite{Gerdts2012}, Theorem 2.3.7), it exists a solution for Problem \ref{LQODEOCP}.
	Moreover, for $\alpha_u , \alpha_p > 0$, the quadratic objective $\hat{J}_{\alpha}$ is even strictly convex due to the included regularization terms and the constraints form a convex set. 
	Hence, the auxiliary problem has a unique solution and, due the convexity of the problem, every local minimum is also a global minimum (see, for instance, \cite{Gerdts2012}, Theorem 2.3.41). \\
\end{remark}

\noindent
Having a solution $(\delta u^{(k)},\delta p^{(k)})$ for the auxiliary Problem \ref{LQODEOCP}, we are able to compute a new iterate $(u^{(k+1)},p^{(k+1)})$ as proposed in \eqref{IterProc}, where the stepsizes $\gamma_k > 0$ in each iteration are computed by an Armijo-like backtracking line search method, see \cite{UlbrichBrüder2011}. 

If we allow for box-constraints in Problem \ref{NonlinearODEOCP}, that is, $\mathcal{U} \subset L^{\infty}(I;\R^{n_u}) \times \R^{n_p}$, in this situation we have to guarantee, that the computed iterates within the Gauss-Newton algorithm are feasible in the sense that $(u^{(k+1)},p^{(k+1)})$ satisfies the box-constraints in each iteration step $k$. We follow the idea of projection methods (see e.g. \cite{UlbrichBrüder2011}) and define the following suitable projection onto $\mathcal{U}$. 

\begin{definition}
(\textit{Projection Operator})\label{DefProjOp} The projection operator that maps each pair $(u,p) \in L^\infty(I;\R^{n_u}) \times \R^{n_p}$ onto the closed convex set $\mathcal{U}$ representing the box-constraints of Problem \ref{NonlinearODEOCP}, is defined as 
\begin{equation}
	\begin{split}
		&\Pi_\mathcal{U}: L^{\infty}(I;\R^{n_u}) \times \R^{n_p} \longrightarrow \mathcal{U} \\	
		&\Pi_\mathcal{U}(u,p)(t) = \Big((\ti{u}_i(t))_{i=1,..,n_u},(\ti{p}_j)_{j=1,..,n_p}\Big),
	\end{split}
\end{equation}
with
\begin{equation}
	\ti u_i(t)=\begin{cases}
		u^\text{low}_i &\mbox{, if } u_i(t) < u^\text{low}_i,\\
		u_i(t)         &\mbox{, if } u^\text{low}_i \leq u_i(t) \leq u^\text{up}_i,\\
		u^\text{up}_i &\mbox{, if } u_i(t) > u^\text{up}_i,
	\end{cases} \quad \text{and} \quad 
	\ti p_j=\begin{cases}
		p^\text{low}_j &\mbox{, if } p_j < p^\text{low}_j,\\
		p_j         &\mbox{, if } p^\text{low}_j \leq p_j \leq p^\text{up}_j,\\
		p^\text{up}_j &\mbox{, if } p_j > p^\text{up}_j.
	\end{cases}
\end{equation} \\
\end{definition} 

\noindent
Bringing together all the concepts discussed, we can now devise a projected Gauss-Newton method, outlined in Algorithm \ref{Gauss-Newton for ODE Optimal Control}, aimed at addressing Problem \ref{NonlinearODEOCP}.

\begin{algorithm}[h]
	\caption{Projected Gauss-Newton Method for ODE Optimal Tracking Control}\label{Gauss-Newton for ODE Optimal Control}
	\begin{algorithmic}
		\State Choose an admissible starting point $(u^{(0)},p^{(0)}) \in \mathcal{U}$, $J_{TOL} > 0$, $\beta,\sigma \in (0,1)$ and set $k=0$.
		\begin{itemize}
			\item[(1)] Solve the dynamical system \eqref{DynSys},\eqref{OutputEqu} for $y^{(k)} = \S(u^{(k)},p^{(k)})$ forwards in time. \newline
			If $\left| J(u^{(k)},p^{(k)}) \right| \leq J_{TOL}$, then STOP. Otherwise go to step (2).
			\item [(2)] Solve the auxiliary problem \ref{LQODEOCP} for $(\delta u^{(k)}, \delta p^{(k)})$.
			\item[(3)] Finally, set $\left(u^{(k+1)},p^{(k+1)}\right) = \Pi_\mathcal{U} \left[ \left(u^{(k)},p^{(k)}\right) + \gamma_k \left(\delta u^{(k)}, \delta p^{(k)}\right) \right]$ whereas the step size $\gamma_k$ is determined by an Armijo-like backtracking line search method, i.e., find the smallest number $i \in \N$ such that
			\begin{align}\label{ArmijoGN}
				J\left( \Pi_\mathcal{U} \left[ (u^{(k)},p^{(k)}) + \beta^{i} (\delta u^{(k)}, \delta p^{(k)}) \right]\right) \leq J(u^{(k)},p^{(k)}) + \sigma \beta^{i} J^\prime (u^{(k)},p^{(k)}) (\delta u^{(k)}, \delta p^{(k)}) 
			\end{align} 
			\item[(4)] Set $k = k+1$ and go to step (1).
		\end{itemize}
	\end{algorithmic}
\end{algorithm}

\subsection{Solving the auxiliary problems}
In the following, we present two different approaches to solve the auxiliary problems arising in the Gauss-Newton method, Algorithm \ref{Gauss-Newton for ODE Optimal Control}, step (2). First, a gradient-based method is presented that solves the auxiliary problems \ref{LQODEOCP} numerically in each iteration. Second, for a known and fixed parameter vector $p = \bar{p} \in \R^{n_p}$, where we only optimize for a control function $u \in L^{\infty}(I;\R^{n_u})$, an explicit solution formula for the auxiliary problems is given via the Riccati-approach for linear-quadratic ODE optimal control problems. 

\subsubsection{Gradient descent method in function spaces}
Again, the idea is to treat Problem \ref{LQODEOCP} as an infinite-dimensional optimization problem 
and apply a gradient descent method in $L^\infty(I;\R^{n_u}) \times \R^{n_p}$, cf. \cite{Gerdts2012}~. This results in an iterative procedure of the form
\begin{align*}
	\left(\delta u^{(k)}_{i+1},\delta p^{(k)}_{i+1}\right) = \left(\delta u^{(k)}_i, \delta p^{(k)}_i\right) + \gamma_i \cdot d_i, \quad i = 0,1,2,...,
\end{align*}
in which we follow the direction of steepest descent $d_i$ of the cost function $\hat{J}_\alpha$ at the current iterate $(\delta u^{(k)}_i, \delta p^{(k)}_i)$ to compute an update. In the gradient descent method, the descent direction $d_i$ is given by the gradient of $\hat{J}_\alpha$ w.r.t. $(\delta u, \delta p)$.

For the derivation of an expression for the gradient $\nabla_{(\delta u,\delta p)} \hat{J}_\alpha (\delta u^{(k)}, \delta p^{(k)}) \in L^{\infty}(I;\R^{n_u}) \times \R^{n_p}$, we proceed as in \cite{Gerdts2012}. Therefore, we first rewrite the cost function $\hat{J}_\alpha$ of the auxiliary problem in the form
\begin{align}\label{Form_J_alpha}
	\hat{J}_\alpha(\delta u,\delta p) = \hat{\phi}(t_f,\delta x(t_f),\delta u(t_f),\delta p) + \int_{t_0}^{t_f} \hat{\varphi}(t,\delta x(t), \delta u(t), \delta p) ~ dt
\end{align}
with 
{\small
\begin{equation}\label{hat_phi_varphi}
	\begin{split}
		& \hat{\phi}(t_f,\delta x(t_f),\delta u(t_f),\delta p) = \frac{1}{2} \left\Vert T^{\frac{1}{2}}\left[ C(t_f)\delta x(t_f) +  D_u(t_f)\delta u(t_f) + D_p(t_f)\delta p - r(t_f) \right] \right\Vert_{\R^{n_y}}^2, \\
		& \hat{\varphi}(t,\delta x,\delta u,\delta p) = \frac{1}{2} \left\Vert Q^{\frac{1}{2}}\left[ C(t)\delta x +  D_u(t)\delta u + D_p(t)\delta p - r(t) \right] \right\Vert_{\R^{n_y}}^2 + \frac{\alpha_u}{2} \left[ 2 {u^{(k)}(t)}^\top \delta u + \delta u^\top \delta u \right] + \frac{\alpha_p}{2(t_f-t_0)} \left[ 2 {p^{(k)}}^\top \delta p + \delta p^\top \delta p \right],
	\end{split}
\end{equation}}
similar to $J$ in Section \ref{sec23}. Then, we consider the corresponding Hamilton function $\mathcal{H}: I \times \R^{n_x} \times \R^{n_u} \times \R^{n_p} \times \R^{n_x} \longrightarrow \R^{n_x}$ that is given by
\begin{align}\label{HamiltonFunction}
	\mathcal{H}(t,\delta x, \delta u, \delta p, \lambda) := ~ \hat{\varphi}(t,\delta x,\delta u,\delta p) + \lambda^\top \left[ A(t) \delta x + B_u(t) \delta u + B_p(t) \delta p \right].
\end{align}

Using the Hamilton function $\mathcal{H}$ and the corresponding adjoint equations \eqref{AdjointEquations} for Problem \ref{LQODEOCP} resulting from the necessary optimality conditions, we get an expression for $\nabla_{(\delta u,\delta p)} \hat{J}_\alpha (\delta u^{(k)}, \delta p^{(k)})$ as stated in the following definition.
\begin{definition}
	Let $(\hat{\delta u}, \hat{\delta p}) \in L^{\infty}(I;\R^{n_u}) \times \R^{n_p}$ be given with $\hat{\delta x} = \ti{\S}^\prime (u^{(k)},p^{(k)})(\hat{\delta u}, \hat{\delta p})$ 
	and let $\lambda \in W^{1,\infty}(I;\R^{n_x})$ satisfy the adjoint equations
	\begin{align}\label{AdjointEquations}
		\begin{cases}
			\dot{\lambda}(t) &= - \mathcal{H}_x(t,\hat{\delta x}(t), \hat{\delta u}(t), \hat{\delta p},\lambda (t))^\top ~~ t \in I, \\
			\lambda(t_f) &= \hat{\phi}_x [t_f]^\top,          
		\end{cases}
	\end{align}
	with $\hat{\phi}_x [t_f] = \hat{\phi}_x (t_f, \hat{\delta x}(t_f), \hat{\delta u}(t_f), \hat{\delta p})$.
	The gradient $\nabla_{(\delta u,\delta p)} \hat{J}_\alpha (\hat{\delta u}, \hat{\delta p}) \in L^{\infty}(I;\R^{n_u}) \times \R^{n_p}$ is then defined by
	\begin{align}\label{Gradient}
		&\left[ \nabla_{(\delta u,\delta p)} \hat{J}_\alpha (\hat{\delta u}, \hat{\delta p}) \right] (t) = \left(\begin{array}{c} \mathcal{H}_u(t,\hat{\delta x}(t), \hat{\delta u}(t), \hat{\delta p},\lambda (t))^\top \\ \int_{t_0}^{t_f} \mathcal{H}_p(s,\hat{\delta x}(s), \hat{\delta u}(s), \hat{\delta p},\lambda (s))^\top ~ds \end{array}\right). 
	\end{align}
\end{definition} 

\noindent
The gradient descent method for Problem \ref{LQODEOCP} is summarized in Algorithm \ref{GDMethod}. Similar to the presented Gauss-Newton method, we combine the method with an Armijo-like backtracking line search method for the determination of the stepsize $\gamma_i$  in iteration $i$. 
For a more detailed derivation of a gradient descent method for general ODE optimal control problems and further details also about the convergence of the method we again refer to \cite{Gerdts2012}.

\begin{algorithm}[h]
	\caption{Gradient Descent Method for the Auxiliary Problems}\label{GDMethod}
	\begin{algorithmic}
		\State Choose an admissible starting point $(\delta u^{(k)}_i, \delta p^{(k)}_i) \in L^{\infty}(I;\R^{n_u}) \times \R^{n_p}$, $\hat{J}_{TOL} > 0$, $\hat{\beta}, \hat{\sigma} \in (0,1)$ and set $i=0$.
		\begin{itemize}
			\item[(a)] Solve the following ODE-IVP, given by the constraints \eqref{AuxiliaryProb_ODEConstraints} in Problem \ref{LQODEOCP}~, for $\delta x_i$ forwards in time: 
			\begin{align}
				\begin{cases}
					\dot{\delta x_i}(t) &= A(t)\delta x_i(t) + B_u(t)\delta u^{(k)}_i(t) + B_p(t) \delta p^{(k)}_i ~~ t \in I, \\
					\delta x_i(t_0)&= 0.         
				\end{cases}
			\end{align}
			
			If $\left| \hat{J}_\alpha(\delta u^{(k)}_i, \delta p^{(k)}_i) \right| \leq \hat{J}_{TOL}$, then STOP. Otherwise go to step (b).
			\item [(b)] Compute the gradient $\nabla_{(\delta u,\delta p)} \hat{J}_\alpha (\delta u^{(k)}_i, \delta p^{(k)}_i)$ by use of the Hamilton function $\mathcal{H}$ \eqref{HamiltonFunction}: First solve the adjoint equations \eqref{AdjointEquations} for $\lambda_i$ backwards in time; second set the search direction $d_i$ as given in \eqref{Gradient}:
			\begin{align*}
				d_i := -\nabla_{(\delta u,\delta p)} \hat{J}_\alpha (\delta u^{(k)}_i, \delta p^{(k)}_i) 
			\end{align*}
			\item[(3)] Finally, set $\left(\delta u^{(k)}_{i+1},\delta p^{(k)}_{i+1}\right) = \left(\delta u^{(k)}_i, \delta p^{(k)}_i\right) + \gamma_i d_i$ whereas the step size $\gamma_i$ is determined by an Armijo-like backtracking line search method, i.e., find the smallest number $m \in \N$ such that
			\begin{align}\label{ArmijoGD}
				\hat{J}_\alpha \left(  (\delta u^{(k)}_i, \delta p^{(k)}_i) + \hat{\beta}^{m} d_i \right) \leq \hat{J}_\alpha(\delta u^{(k)}_i, \delta p^{(k)}_i) + \hat{\sigma} \hat{\beta}^{m} \hat{J}_\alpha^\prime (\delta u^{(k)}_i, \delta p^{(k)}_i) (d_i) 
			\end{align} 
			\item[(4)] Set $i = i+1$ and go to step (a).
		\end{itemize}
	\end{algorithmic}
\end{algorithm} 

\begin{remark}
	To tackle the auxiliary problem, presented as a linear-quadratic ODE optimal control challenge in Algorithm \ref{Gauss-Newton for ODE Optimal Control}, step (2), various solution techniques from ODE optimal control literature can be employed, as elaborated in \cite{Gerdts2012}. One approach involves first discretizing the optimal control problem, followed by applying methods from finite-dimensional constrained optimization. However, this method may encounter challenges, especially when dealing with fine time grids or large time intervals, leading to optimization problems with numerous variables that can be computationally demanding (known as the curse of dimensionality).
	Alternatively, an indirect solution approach, falling under the category of first-optimize-then-discretize methods, leverages necessary optimality conditions derived from the corresponding Hamilton function. Yet, numerically solving the resulting system of equations for a solution of the auxiliary problem typically involves addressing a mixed boundary value problem with coupled equations, presenting its own set of difficulties.
	Thus, we decide for a function space method. First-order descent methods, representing the simplest approach in this category, merely necessitates the first derivative for an iterative solution process. Higher-order methods, such as the Newton method, offer potentially faster convergence rates but are impractical due to their requirement for higher derivatives, resulting in more complex operators that may not be readily determined or beneficial for our specific applications.
\end{remark}

\subsection{Riccati-based method}
In cases where we fix the parameter vector $p = \bar{p} \in \R^{n_p}$, Problem \ref{ODEOptControlProblem} reduces to finding a control function $u \in \mathcal{V} \subseteq L^{\infty}(I;\R^{n_u})$,
\begin{align}\label{ConvexSetConstraints_nur_u}
	\mathcal{V} := \left\lbrace u \in L^{\infty}(I;\R^{n_u}) ~\big|~ u^\text{low}_i\leq u_i(t)\leq u^\text{up}_i \mbox{ a.e. in }I, ~i=1,..,n_u \right\rbrace,
\end{align} 
such that the error
\begin{align}\label{ODEOptControl_u}
	J(u) &= \frac{1}{2} \left\Vert Q^{\frac{1}{2}} \left[\S(u,\bar{p}) - y_\text{ref}\right] \right\Vert_{L^2(I;\R^{n_y})}^2 + \frac{\alpha_u}{2} \left\Vert u \right\Vert_{L^2(I;\R^{n_u})}^2  + \frac{1}{2} \left\Vert T^{\frac{1}{2}} \left[\S(u,\bar{p})(t_f) - y_\text{ref}(t_f)\right] \right\Vert_{R^{n_y}}^2 
\end{align}
is minimized. 
Applying the Gauss-Newton approach to the reduced problem again ends with an iterative procedure, where the subproblems \ref{LQODEOCP} in each iteration reduce to finding $\delta u^{(k)} \in L^{\infty}(I;\R^{n_u})$ that solves the following linear-quadratic ODE optimal control problem where all appearing derivatives in the resulting equations with respect to $p$ are equal zero: 
{\small
\begin{equation}\label{RiccatiAuxiliaryProblem_u}
	\begin{split}
		\underset{\delta u^{(k)} \in L^{\infty}(I;\R^{n_u})}{\min} ~~  \hat{J}_{\alpha}(\delta u^{(k)}) = & ~~ \int_{t_0}^{t_f} \frac{1}{2} \left\Vert Q^{\frac{1}{2}}\left[ C(t)\delta x(t) +  D_u(t)\delta u^{(k)}(t) - r(t) \right] \right\Vert_{\R^{n_y}}^2 + \frac{\alpha_u}{2} \left[ 2 u^{(k)}(t)^\top \delta u^{(k)}(t) + \left\Vert \delta u^{(k)}(t)\right\Vert_{\R^{n_u}}^2 \right] ~dt\\
		& +  \frac{1}{2} \left\Vert T^{\frac{1}{2}}\left[C(t_f)\delta x(t_f) +  D_u(t_f)\delta u^{(k)}(t_f) - r(t_f) \right] \right\Vert^2_{\R^{n_y}}  \\
		\\		
		\mbox{s.t.}  ~~
		&  \begin{cases}
			\dot{\delta x}(t) &= A(t)\delta x(t) + B_u(t)\delta u^{(k)}(t) ~~ t \in I, \\
			\delta x(t_0)&= 0.         
		\end{cases}	
	\end{split}
\end{equation}}
with matrix functions 
\begin{align*}
	A(t) &= f_x(t,x^{(k)}(t),u^{(k)}(t),\bar{p}) \in \R^{n_x \times n_x}, \quad B_u(t) = f_u(t,x^{(k)}(t),u^{(k)}(t),\bar{p}) \in \R^{n_x \times n_u}, \\
	C(t) &= h_x(t,x^{(k)}(t),u^{(k)}(t),\bar{p}) \in \R^{n_y \times n_x}, \quad D_u(t) = h_u(t,x^{(k)}(t),u^{(k)}(t),\bar{p}) \in \R^{n_y \times n_u}, 
\end{align*} 
$x^{(k)} = \tilde{\S}(u^{(k)},\bar{p})$ and  $r(t) := y_\text{ref}(t) - \S(u^{(k)},\bar{p})(t) = y_\text{ref}(t) - y^{(k)}(t)$. \\

Instead of applying the gradient descent method, Algorithm \ref{GDMethod}, to solve the auxiliary problem \eqref{RiccatiAuxiliaryProblem_u} iteratively, the solution in the case of fixed parameters 
can be derived explicitly by applying the Riccati-approach for linear-quadratic ODE optimal control problems. This approach gives a unique solution $\delta u^{(k)}$ for the problem under some assumptions and, in addition, an explicit solution formula for $\delta u^{(k)}$ in form of a feedback control law including the state $\delta x$ and a matrix-valued function solving the corresponding Riccati differential equations. For more details on the Riccati-approach we refer the reader to \cite{Locatelli2001}, \cite{Sontag1998}. Before applying the results from the Riccati theorem, problem \eqref{RiccatiAuxiliaryProblem_u} first has to be rewritten using appropriate transformations, see, e.g., \cite{Locatelli2001}. We end up with a problem of the form

\begin{equation}\label{TransformedRiccatiAuxiliaryProblem_u}
	\begin{split}
		\underset{\delta v \in L^{\infty}(I;\R^{n_u})}{\min} & ~~ \widehat{\delta x}(t_f)^\top \hat{T} \widehat{\delta x}(t_f) + \int_{t_0}^{t_f} \widehat{\delta x} (t)^\top \hat{Q}(t) \widehat{\delta x}(t) + \delta v(t)^\top R(t) \delta v(t) ~dt\\		
		\mbox{s.t.}  ~~
		&  \begin{cases}
			\dot{\widehat{\delta x}}(t) &= \hat{A}(t)\widehat{\delta x}(t) + \hat{B_u}(t)\delta v(t) ~~ t \in I, \\
			\widehat{\delta x}(t_0) &= \left(\begin{array}{c}
				0 \\ 1
			\end{array}\right),         
		\end{cases}	
	\end{split}
\end{equation}

\noindent
where the resulting weight-matrices for the state $\widehat{\delta x}$ and the control $\delta v$ in the objective function of the transformed problem \eqref{TransformedRiccatiAuxiliaryProblem_u} are of the form
\begin{align}
	\hat{Q}(t) &= \left(\begin{array}{c c}
		C(t)^\top \tilde{Q}(t) C(t) & -C(t)^\top \tilde{Q}(t) r(t) - C(t)^\top R(t)^{-1} Z(t)^\top k(t) \\
		-r(t)^\top \tilde{Q}(t)C(t)  -k(t)^\top Z(t)^\top R(t)^{-1} C(t) & ~~~ r(t)^\top \tilde{Q}(t)r(t) + 2r(t)^\top Z(t) R(t)^{-1}k(t) - k(t)^\top R(t)^{-1} k(t)
	\end{array}\right), \label{HatQ} \\
	\notag \\
	\hat{T} &= \left(\begin{array}{c c}
		C(t_f)^\top \tilde{T} C(t_f) & -C(t_f)^\top \tilde{T} r(t_f)  \\
		-r(t_f)^\top \tilde{T} C(t_f)  & ~~~ r(t_f)^\top \tilde{T} r(t_f) 
	\end{array}\right), \label{HatT}\\
	\notag \\
	R(t) &= D_u(t)^\top Q D_u(t) + \alpha_u \mathbf{1}_{n_u \times n_u} \label{R}
\end{align}
with  
\begin{align*}
	\tilde{Q}(t) &= Q - Z(t) R(t) Z(t)^\top, \quad Z(t) = Q D_u(t), \quad k(t) = \alpha_u u^{(k)}(t), \\
	\tilde{T} &= T - \tilde{Z}(t_f)\tilde{R}(t_f)^{-1} \tilde{Z}(t_f)^\top, \quad \tilde{R}(t_f) = D_u(t_f)^\top T D_u(t_f), \quad \tilde{Z}(t_f) = T D_u(t_f).
\end{align*}

In the following, we have to make assumptions for the matrices $\hat{Q}(t)$ and $\hat{T}$ to satisfy all preconditions for the Riccati theorem before applying the results. 

\begin{assumption}
\label{ass3} The matrix $\hat{Q}(t)$ is assumed to be positive semi-definite for all $t \in I$. Further, $\hat{T}$ is also assumed to be positive semi-definite.
\end{assumption}

\begin{remark}
	Obviously, $\hat{T}$ is symmetric and $\hat{Q}(t)$ and $R(t)$ are symmetric for all $t \in I$. For $\alpha_u > 0$ the required positive definiteness of the matrix $R(t)$ for all $t \in I$ can be ensured while the positive semi-definiteness of the matrices $\hat{Q}(t)$ and $\hat{T}$ is assumed in the following and its verification in a specific application scenario is left to the reader.\\
\end{remark}

\noindent
In the end, by applying Riccati to the transformed problem \eqref{TransformedRiccatiAuxiliaryProblem_u} and performing the corresponding back-transformations, the solution of the auxiliary problem \eqref{RiccatiAuxiliaryProblem_u} is given by the control law
\begin{equation}\label{Optimal Control Law_u}
	\delta u^{(k)}(t) = \left(F(t) - R(t)^{-1}Z(t)^\top C(t)\right)\delta x(t)	- R(t)^{-1}B_u(t)^\top \beta(t) + R(t)^{-1}Z(t)^\top r(t) - R(t)^{-1}k(t)
\end{equation}
with
\begin{align*}
	F(t) = -R(t)^{-1}B_u(t)^\top P(t),
\end{align*}
where the symmetric matrix-function $P(t) \in \R^{n_x \times n_x}$ solves the corresponding Riccati differential equations
\begin{equation}\label{Riccati Equation for Auxiliary Problem_u}
	\begin{split}
		\begin{cases}
			\dot{P}(t) &= -P(t)\left(A(t)-B_u(t)R(t)^{-1}Z(t)^\top C(t) \right) - \left(A(t)-B_u(t)R(t)^{-1}Z(t)^\top C(t) \right)^\top P(t) \\
			& ~~~ + P(t)B_u(t)R(t)^{-1} B_u(t)^\top P(t) - C(t)^\top \tilde{Q}(t) C(t), \\
			P(t_f) &= C(t_f)^\top \tilde{T} C(t_f),
		\end{cases}
	\end{split}
\end{equation}
and $\beta(t) \in \R^{n_x}$ solves the end-value problem
\begin{equation}\label{Beta Equation for Auxiliary Problem_u}
	\begin{split}
		\begin{cases}
			\dot{\beta}(t) &= -\left(A(t) - B_u(t)R(t)^{-1}Z(t)^\top C(t) + B_u(t) F(t)\right)^\top \beta(t) - P(t) \varphi(t) \\
			& ~~~ + C(t)^\top \tilde{Q}(t) r(t) + C(t)^\top Z(t) R(t)^{-1} k(t)\\
			\beta(t_f) &= -C(t_f)^\top \tilde{T} r(t_f),
		\end{cases}
	\end{split}
\end{equation}
with
\begin{align*}
	\varphi(t) = - B_u(t) R(t)^{-1} k(t) + B_u(t) R(t)^{-1} Z(t)^\top r(t),
\end{align*}
both on $I$. Finally, the solution procedure based on the Riccati-approach for the auxiliary problem \eqref{RiccatiAuxiliaryProblem_u}, is summarized in Algorithm \ref{SolutionProcedure}.

\begin{algorithm}[h]
	\caption{Riccati-Based Method for the Auxiliary Problems}\label{SolutionProcedure} 
	\begin{algorithmic}
		\State
		\begin{itemize}
			\item [(a)] Solve the Riccati equations \eqref{Riccati Equation for Auxiliary Problem_u} for the $(n_x \times n_x)$-matrix function $P^{(k)}$ (column-wise) backwards in time. 
			\item [(b)] Solve the equations \eqref{Beta Equation for Auxiliary Problem_u} for the $(n_x \times 1)$-function $\beta^{(k)}$ backwards in time.
			\item[(c)] Define $\delta u^{(k)}$ by the control law \eqref{Optimal Control Law_u} using the results from step (a) and (b), replace $\delta u^{(k)}$ in the linear ODE-constraints of the auxiliary problem \eqref{RiccatiAuxiliaryProblem_u} by the control law and solve the resulting system for $\delta x^{(k)}$ forwards in time:
			\begin{align*}
				\begin{cases}
					\dot{\delta x}^{(k)}(t) &= (A(t) + B_u(t)F(t) - B_u(t) R(t)^{-1}Z(t)^\top C(t) ) \delta x^{(k)}(t) - B_u(t) R(t)^{-1} B_u(t)^\top \beta^{(k)}(t)\\
					&\quad + B_u(t) R(t)^{-1}Z(t)^\top r(t) - B_u(t) R(t)^{-1} k(t),   \\
					\delta x^{(k)}(t_0) &= 0,	
				\end{cases}
			\end{align*}
			\item[(d)] Compute the optimal control $\delta u^{(k)}$ with \eqref{Optimal Control Law_u} using $\delta x^{(k)}$, $P^{(k)}$ and $\beta^{(k)}$ from steps (a)-(c). 
		\end{itemize}
	\end{algorithmic}
\end{algorithm} 

\begin{remark}
	The Riccati theory, employed to address linear quadratic ODE optimal control problems under specific assumptions, falls within the realm of indirect solution methods, as detailed in \cite{Gerdts2012}. This theory is typically presented for controls $\delta u \in L^{\infty}(I,\R^{n_u})$. For constant control functions, such as $\delta p \in \R^{n_p}$, existence and uniqueness statements can still be provided. However, an explicit solution formula exists only under additional assumptions. To obtain such a solution, one must solve not only the Riccati differential equations but also the corresponding sensitivity equations \cite{Gerdts2012} for \eqref{DynSys}, alongside an integro-differential equation for $\delta x$ as discussed in \cite{Lakshmikantham1995}.
	The numerical solution of this coupled system of equations is typically challenging. Furthermore, extending the theory to accommodate problems involving both types of input functions, as required in our scenario with $\delta u$ and $\delta p$, necessitates addressing a system of coupled equations akin to when considering the necessary optimality conditions provided through the Hamilton function.
\end{remark}

\subsection{Convergence}
In this section, we establish a convergence theorem for the Gauss-Newton method under the assumption of fixed parameters $p = \bar{p} \in \R^{n_p}$. Here, the solution of the auxiliary problem is explicitly determined through the necessary optimality conditions. Instead of employing the Riccati-approach, we utilize the normal equations corresponding to the auxiliary problem. It is worth noting that both methodologies are equivalent; the solution obtained via the Riccati-approach also satisfies the normal equations, and vice versa. Furthermore, both approaches stem from the consideration of the first-order necessary optimality conditions for the auxiliary problem.
We further assume for the weight matrices that $Q = I_{n_y \times n_y}$ and $T = 0_{n_y \times n_y}$ and that no box-constraints for the control are included, i.e. $\mathcal{V} = L^{\infty}(I;\R^{n_u})$, such that we do not need a projection operator.

\subsubsection{Adjoint Operator and Normal Equations  }
We have already shown in Subsection \ref{sec23} that $\S(\cdot, \bar{p}) :L^\infty(I;\R^{n_u}) \rightarrow L^\infty(I;\R^{n_y})$ for fixed $\bar{p} \in \R^{n_p}$ is Fr\'echet-differentiable with the Fr\'echet-derivative $\S^\prime(\hat{u},\bar{p}) :L^\infty(I;\R^{n_u}) \rightarrow L^\infty(I;\R^{n_y})$ at $\hat{u} \in L^\infty(I;\R^{n_u})$ being a bounded linear operator and the solution operator of a linear state-space model. 
The adjoint operator, see \cite{LiYong1995}, \cite{UlbrichBrüder2011}, of $\S^\prime(\hat{u},\bar{p})$ is given as follows.
\begin{theorem}\label{Theorem_Adjoint}
	Let $\hat{u} \in L^{\infty}(I;\R^{n_u})$ with $\tilde{\S}(\hat{u},\bar{p}) = \hat{x}$ be given. The adjoint operator of the map $\S^\prime(\hat{u},\bar{p})$ is then defined as
	\begin{equation}
		\left(\S^\prime(\hat{u},\bar{p})\right)^\ast:  L^{\infty}(I; \R^{n_y}) \longrightarrow L^{\infty}(I;\R^{n_u}) 
	\end{equation}
	with $\left(\S^\prime(\hat{u},\bar{p})\right)^\ast(\delta y) = \delta u$ and is the solution operator of the following linear system
	\begin{align}
		\dot{\delta x}(t)  &= - f_x[t]^\top \delta x(t) - h_x[t]^\top \delta y(t), ~~ \delta x (t_f) = 0 \label{AdjDynSys1}\\
		\delta u(t) &= f_u[t]^\top \delta x(t) + h_u[t]^\top \delta y(t) \label{AdjDynSys2}
	\end{align}
	containing an ODE end-value problem for $\delta x \in W^{1,\infty}(I;\R^{n_x})$.
\end{theorem}

\begin{remark}
	Since $I = [t_0,t_f]$ is compact, it holds $x \in L^\infty(I;\R^{n}) \subset L^2(I;\R^{n})$ and we can thus use the corresponding inner product of $L^2$ for the derivation of the adjoint map, see Appendix for the proof.\\
\end{remark}

\noindent
The unique solution of the subproblems is, besides of the Riccati-approach, also given by the solution of the normal equations that may be derived by applying the necessary first-order optimality conditions.  

\begin{lemma}
The solution of Problem \ref{LQODEOCP} for fixed $\bar{p} \in \R^{n_p}$ is given by
\begin{equation}\label{NormalEquations}
	\left[ \S^\prime(u^{(k)},\bar{p})^\ast \S^\prime(u^{(k)},\bar{p}) + \alpha_u \mathbf{I} \right] \left(\delta u^{(k)}\right) = - \left[ \S^\prime(u^{(k)},\bar{p})^\ast \left[\S(u^{(k)},\bar{p}) - y_\text{ref}\right] + \alpha_u u^{(k)} \right] \quad \left( \text{in}~ L^2(I;\R^{n_u}) \right). 
\end{equation}
\end{lemma} 

\noindent
For $\delta u^{(k)} \in L^\infty(I;\R^{n_u})$ one can guarantee that the operator $\left[ \S^\prime(u^{(k)},\bar{p})^\ast \S^\prime(u^{(k)},\bar{p}) + \alpha_u \mathbf{I} \right]$ maps to $L^\infty(I;\R^{n_u})$ and is invertible for $\alpha_u > 0$. Note that the map and its inverse are both solution operators of a state-space representation including a linear ODE boundary value problem for which an existence and uniqueness result for a solution can be found in \cite{Ascher1988}. Thus, although equation \eqref{NormalEquations} holds generally in the $L^2$-sense, we explicitly solve it for $\delta u^{(k)} \in L^{\infty}(I;\R^{n_u})$ as stated in \eqref{SolNormalEquations}. 

\subsubsection{Convergence Result}
The objective is to demonstrate that the method outlined is a descent method within the space $L^\infty(I;\R^{n_u})$, with convergence guaranteed in the space $L^2(I;\R^{n_u})$. To establish this, we rely on a global convergence theorem (Theorem 2.2, \cite{UlbrichBrüder2011}) for admissible directions and stepsizes.

\begin{theorem}  
Let, in addition to the Assumptions \ref{Ass1}, \ref{Ass2}, the functions $f$ and $h$ to be chosen such that $J^\prime$ is uniformly continuous on the set $N_0^\rho := \{ u + \delta u~ |~ J(u) \leq J(u^{(0)}), \Vert\delta u\Vert_{L^2(I;\R^{n_u})} \leq \rho\}$ for some $\rho > 0$. Furthermore, consider the sequence $\{u^{(k)}\}$ generated by the Gauss-Newton method, Algorithm \ref{Gauss-Newton for ODE Optimal Control}, with $\{\gamma_k \}$ being the sequence of Armijo stepsizes. Let, in addition, the sequence $\{\delta u^{(k)}\}$ be given by the solutions of the normal equations \eqref{NormalEquations}, i.e., 
\begin{align}\label{SolNormalEquations}
	\delta u^{(k)} = - \left[ \S^\prime(u^{(k)},\bar{p})^\ast \S^\prime(u^{(k)},\bar{p}) + \alpha_u \mathbf{I} \right]^{-1} \left[ \S^\prime(u^{(k)},\bar{p})^\ast \left[\S(u^{(k)},\bar{p}) - y_\text{ref}\right] + \alpha_u u^{(k)} \right] \in L^\infty(I;\R^{n_u}) 
\end{align} 
and such that they are not too short in the following sense:
\begin{align}\label{LowBoundLengthStep2} 
	\Vert \delta u^{(k)} \Vert_{L^2(I;\R^{n_u})} \geq - \frac{J^\prime (u^{(k)})(\delta u^{(k)})}{	\Vert \delta u^{(k)} \Vert_{L^2(I;\R^{n_u})}}.
\end{align} 
Then, the sequences $\{\delta u^{(k)}\}$ and $\{u^{(k)}\}$ have the following properties:
\begin{itemize}
	\item[(1)] For all $k$ it holds that $\delta u^{(k)}$ is a descent direction, i.e., $J^\prime(u^{(k)})(\delta u^{(k)}) < 0$.
	\item[(2)] The sequence of function values $\{J(u^{(k)})\}$ is monotonically decreasing, i.e., for all $k$ it holds $J(u^{(k+1)}) < J(u^{(k)})$.
	\item[(3)] Every accumulation point of $\{u^{(k)}\}$ is a stationary point of $J$. 
\end{itemize}
\end{theorem}

\begin{proof}
For the proof we verify that all preconditions from Theorem 2.2 stated in \cite{UlbrichBrüder2011} are satisfied for the considered setting. \\
Under Assumptions \ref{Ass1} and \ref{Ass2}, as already shown in Section \ref{sec23}, the objective function $J: L^{\infty}(I;\R^{n_u}) \longrightarrow \R$ is (continuously) Fr\'echet-differentiable and its Fr\'echet-derivative  $J^\prime (u^{(k)}): L^{\infty}(I;\R^{n_u}) \longrightarrow \R$ at arbitrary $u^{(k)} \in L^{\infty}(I;\R^{n_u})$ is, by definition, a bounded linear operator and can also be written as follows using the Fr\'echet-derivative of $\S$ and its adjoint:
\begin{align}\label{FréchetDerv}
	J^\prime(u^{(k)})(\delta u) = \langle \S^\prime(u^{(k)},\bar{p})^\ast \left[\S(u^{(k)},\bar{p}) - y_\text{ref}\right] + \alpha_u u^{(k)}, \delta u \rangle_{L^{2}(I;\R^{n_u})}.
\end{align}

\noindent
\textit{Proof of (1):} The operator $\S^\prime(u^{(k)},\bar{p})^\ast \S^\prime(u^{(k)},\bar{p}) + \alpha_u \mathbf{I}$ satisfies for arbitrary $\delta u \in L^{\infty}(I;\R^{n_u})$, $\delta u \neq 0$:
\begin{equation}\label{OperProp}
	\begin{split}
		\langle \left[ \S^\prime(u^{(k)},\bar{p})^\ast \S^\prime(u^{(k)},\bar{p}) + \alpha_u \mathbf{I} \right] (\delta u), \delta u \rangle_{L^{2}(I;\R^{n_u})} 
		= &\langle \left[ \S^\prime(u^{(k)},\bar{p})^\ast \S^\prime(u^{(k)},\bar{p}) \right] (\delta u), \delta u \rangle_{L^{2}(I;\R^{n_u})} + \alpha_u \langle \delta u, \delta u \rangle_{L^{2}(I;\R^{n_u})} \\
		= & \Vert \S^\prime(u^{(k)},\bar{p})(\delta u) \Vert_{L^2(I;\R^{n_y})}^2 + \alpha_u \Vert \delta u \Vert_{L^2(I;\R^{n_u})}^2 \\
		\geq & \alpha_u \Vert \delta u \Vert_{L^2(I;\R^{n_u})}^2 \\
		> & 0.
	\end{split}
\end{equation} 
Similarly, for the inverse operator $\left[ \S^\prime(u^{(k)},\bar{p})^\ast \S^\prime(u^{(k)},\bar{p}) + \alpha_u \mathbf{I} \right]^{-1}$ we have 
\begin{equation}\label{InvOperProp}
	\langle \left[ \S^\prime(u^{(k)},\bar{p})^\ast \S^\prime(u^{(k)},\bar{p}) + \alpha_u \mathbf{I} \right]^{-1} (\delta v), \delta v \rangle_{L^{2}(I;\R^{n_u})} > 0
\end{equation}
Thus, we get 
{\small
\begin{equation}\label{DescDirec}
	\begin{split}
		& J^\prime(u^{(k)})(\delta u^{(k)}) \\
		&= \langle \S^\prime(u^{(k)},\bar{p})^\ast \left[\S(u^{(k)},\bar{p}) - y_\text{ref}\right] + \alpha_u u^{(k)}, - \left[ \S^\prime(u^{(k)},\bar{p})^\ast \S^\prime(u^{(k)},\bar{p}) + \alpha_u \mathbf{I} \right]^{-1} \left[ \S^\prime(u^{(k)},\bar{p})^\ast \left[\S(u^{(k)},\bar{p}) - y_\text{ref}\right] + \alpha_u u^{(k)} \right] \rangle_{L^{2}(I;\R^{n_u})} \\
		&= - \langle \ti u, \left[ \S^\prime(u^{(k)},\bar{p})^\ast \S^\prime(u^{(k)},\bar{p}) + \alpha_u \mathbf{I} \right]^{-1} \ti u \rangle_{L^{2}(I;\R^{n_u})} \\
		&\overset{\eqref{InvOperProp}}{<} 0
	\end{split}
\end{equation}}
for $\ti u := \S^\prime(u^{(k)},\bar{p})^\ast \left[\S(u^{(k)},\bar{p}) - y_\text{ref}\right] + \alpha_u u^{(k)} \neq 0$. This shows that $\delta u^{(k)}$ is a descent direction for every $k$.\\ 

\noindent
\textit{Proof of (2):} Follows directly from $\delta u^{(k)}$ being a descent direction by considering the Armijo stepsize rule: \begin{equation}\label{DeacrFuncVal}
	\begin{split}
		J(u^{(k)} + \gamma_k \delta u^{(k)}) - J(u^{(k)}) \leq \sigma \gamma_k J^\prime (u^{(k)})(\delta u^{(k)}) \overset{\eqref{DescDirec}}{<} 0 \quad \text{implies} \quad  J(u^{(k+1)}) < J(u^{(k)}) ~\forall k
	\end{split}
\end{equation}
for $\gamma_k = \beta^i$, $\sigma, \beta > 0$ and $u^{(k+1)} = u^{(k)} + \gamma_k \delta u^{(k)}$. \\

\noindent
\textit{Proof of (3):} The continuously Fr\'echet-differentiable objective function $J : L^{\infty}(I;\R^{n_u}) \longrightarrow \R$ is, by definition, bounded below by zero, i.e., $J(u) \geq 0$ for all $u \in L^\infty(I;\R^{n_u})$. We first show the admissibility of the descent directions $\delta u^{(k)}$, that is,   
\begin{align*}
	\frac{- J^\prime(u^{(k)})(\delta u^{(k)})}{\Vert \delta u^{(k)} \Vert_{L^2(I;\R^{n_u})}} \overset{k \rightarrow \infty}{\longrightarrow } 0 \quad \text{implies} \quad  \Vert J^\prime (u^{(k)}) \Vert_{\text{Op}} \overset{k \rightarrow \infty}{\longrightarrow } 0.
\end{align*} 
For the numerator $J^\prime (u^{(k)})(\delta u^{(k)})$, we get with \eqref{SolNormalEquations} and \eqref{FréchetDerv} that
\begin{equation}\label{DervAbsch}
	\begin{split}
		J^\prime (u^{(k)})(\delta u^{(k)}) &= \langle - \left[ \S^\prime(u^{(k)},\bar{p})^\ast \S^\prime(u^{(k)},\bar{p}) + \alpha_u \mathbf{I} \right] \left(\delta u^{(k)}\right) , \delta u^{(k)} \rangle_{L^{2}(I;\R^{n_u})} \\
		& \overset{\eqref{OperProp}}{\leq} - \alpha_u \Vert \delta u^{(k)} \Vert_{L^2(I;\R^{n_u})}^2	
	\end{split}
\end{equation}
for $\delta u^{(k)} \neq 0$. Further, we have for the operator norm of $J^\prime(u^{(k)})$, given by
\begin{equation}
	\begin{split}
		\Vert J^\prime (u^{(k)}) \Vert_{\text{Op}} := \underset{\Vert \delta u \Vert_{L^\infty(I;\R^{n_u})} = 1}{\sup}	\lvert J^\prime (u^{(k)})(\delta u) \rvert,
	\end{split}
\end{equation}
that
\begin{equation}\label{Absch}
	\begin{split}
		\lvert J^\prime (u^{(k)})(\delta u) \rvert^2 & = \bigl| \langle \S^\prime(u^{(k)},\bar{p})^\ast \left[\S(u^{(k)},\bar{p}) - y_\text{ref}\right] + \alpha_u u^{(k)}, \delta u \rangle_{L^{2}(I;\R^{n_u})} \bigr|^2 \\
		& \overset{\text{C.S.}}{\leq} \Vert \S^\prime(u^{(k)},\bar{p})^\ast \left[\S(u^{(k)},\bar{p}) - y_\text{ref}\right] + \alpha_u u^{(k)} \Vert_{L^2(I;\R^{n_u})}^2 \cdot \Vert \delta u \Vert_{L^2(I;\R^{n_u})}^2 \\
		& \leq \Vert \S^\prime(u^{(k)},\bar{p})^\ast \left[\S(u^{(k)},\bar{p}) - y_\text{ref}\right] + \alpha_u u^{(k)} \Vert_{L^2(I;\R^{n_u})}^2 \cdot \Vert \delta u \Vert_{L^\infty(I;\R^{n_u})}^2,
	\end{split}
\end{equation}
where we used the Cauchy-Schwarz (C.S.) inequality. Hence, 
\begin{equation}\label{Absch2}
	\Vert J^\prime (u^{(k)}) \Vert_{\text{Op}} := \underset{\Vert \delta u \Vert_{L^\infty(I;\R^{n_u})} = 1}{\sup}	\lvert J^\prime (u^{(k)})(\delta u) \rvert \overset{\eqref{Absch}}{\leq} \Vert \S^\prime(u^{(k)},\bar{p})^\ast \left[\S(u^{(k)},\bar{p}) - y_\text{ref}\right] + \alpha_u u^{(k)} \Vert_{L^2(I;\R^{n_u})}
\end{equation}
Using again \eqref{SolNormalEquations}, we further follow:
\begin{equation}\label{OpNorm}
	\begin{split}
		\Vert J^\prime (u^{(k)}) \Vert_{\text{Op}} & \overset{\eqref{Absch2}}{\leq} \Vert \S^\prime(u^{(k)},\bar{p})^\ast \left[\S(u^{(k)},\bar{p}) - y_\text{ref}\right] + \alpha_u u^{(k)} \Vert_{L^2(I;\R^{n_u})} \\
		& \leq \Vert - \left[ \S^\prime(u^{(k)},\bar{p})^\ast \S^\prime(u^{(k)},\bar{p}) + \alpha_u \mathbf{I} \right] \left(\delta u^{(k)} \right) \Vert_{L^2(I;\R^{n_u})} \\
		& \leq \Vert \S^\prime(u^{(k)},\bar{p})^\ast \S^\prime(u^{(k)},\bar{p}) + \alpha_u \mathbf{I} \Vert_{\text{Op}} \cdot \Vert \delta u^{(k)} \Vert_{L^2(I;\R^{n_u})} \\
		& \leq M \cdot \Vert \delta u^{(k)} \Vert_{L^2(I;\R^{n_u})} \\
	\end{split}
\end{equation}
where we used that $\S^\prime(u^{(k)},\bar{p})^\ast \S^\prime(u^{(k)},\bar{p}) + \alpha_u \mathbf{I}$ is a bounded operator, i.e., $\Vert \S^\prime(u^{(k)},\bar{p})^\ast \S^\prime(u^{(k)},\bar{p}) + \alpha_u \mathbf{I} \Vert_{\text{Op}} \leq M$. \\
Finally, using \eqref{DervAbsch} and \eqref{OpNorm}, we have the following inequality chain:
\begin{equation}
	\frac{- J^\prime(u^{(k)})(\delta u^{(k)})}{\Vert \delta u^{(k)} \Vert_{L^2(I;\R^{n_u})}} \overset{\eqref{DervAbsch}}{\geq} \frac{\alpha_u \Vert \delta u^{(k)} \Vert_{L^2(I;\R^{n_u})}^2}{\Vert \delta u^{(k)} \Vert_{L^2(I;\R^{n_u})}} = \alpha_u \Vert \delta u^{(k)} \Vert_{L^2(I;\R^{n_u})}^2 \overset{\eqref{OpNorm}}{\geq} \frac{\alpha_u}{M} \cdot \Vert J^\prime (u^{(k)}) \Vert_{\text{Op}} > 0
\end{equation} 
If now, $\frac{- J^\prime(u^{(k)})(\delta u^{(k)})}{\Vert \delta u^{(k)} \Vert_{L^2(I;\R^{n_u})}} \longrightarrow 0$ for $k \rightarrow \infty$, it follows that $\Vert J^\prime (u^{(k)}) \Vert_{\text{Op}} \overset{k \rightarrow \infty}{\longrightarrow } 0$ and, thus, the admissibility of the search directions. \\
For the proof of the admissibiliy of the Armijo stepsizes we refer to \cite{UlbrichBrüder2011}, p. 101ff.
Finally, with the admissibility of the descent directions and stepsizes, for the proof of (3) one can follow the proof of Theorem 2.2, \cite{UlbrichBrüder2011}~, which shows the desired statement.\\
\end{proof}

\begin{remark}
	For the proof we considered the convergence of the sequence $\{u^{(k)}\}$ with respect to the $L^2$-norm. The algorithm starts with $u^{(0)} \in L^{\infty}(I;\R^{n_u})$ and, since $\delta u^{(k)} \in L^\infty(I;\R^{n_u})$, all further iterates $u^{(k)}$ stay in $L^\infty(I;\R^{n_u})$. Under the assumption that the iterates $u^{(k)}$ are uniformly bounded in $L^\infty(I;\R^{n_u})$, every accumulation point is also in $L^\infty(I;\R^{n_u})$. 
\end{remark}

\section{Simulating Example - Nonlinear Quarter-Car Model with Acceleration of the Upper Body Mass as Output}\label{sec4}
\begin{figure}[h]
	\centering
	\includegraphics[scale=0.7]{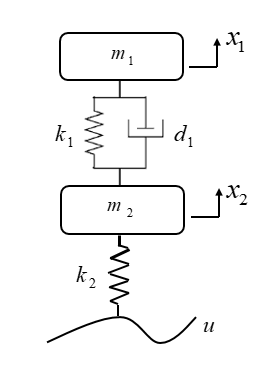}
	\caption{Structure of the quarter-car model with parameter notation.}\label{Figure QCM}
\end{figure}

We consider as an example for a dynamical system $\S$ a nonlinear quarter-car model, see Figure \ref{Figure QCM}, given by a one-dimensional mass-spring-damper system with two bodies. The upper body typically describes a quarter of a vehicle's chassis while the lower body depicts the wheel and the rim. These bodies are connected via spring-damper elements that, in our example, are assumed to have a nonlinear force-law for the spring and a linear one for the damper. The lower body is connected to the ground by a linear spring. 
We simulate a drive of the quarter-car over a road profile at constant speed, that is,  the road profile becomes a time-dependent input function $u: I \rightarrow \R$, $I \subset \R$, $n_u = 1$ and, additionally in our setting, we consider an unknown model parameter $p = k_1 \in \R$, $n_p = 1$ representing the stiffness of the spring between the upper and lower body mass. 
We denote with $x_1, x_2: I \rightarrow \R$ and $x_3, x_4: I \rightarrow \R$ the displacements and the velocities of the upper and lower bodies, i.e., $n_x = 4$, and we consider as an output quantity of interest the time-dependent acceleration of the upper body mass referred to as $y: I \rightarrow \R $, $n_y = 1$.

\noindent
The model equations for the quarter-car model derived via Newton's second law, see, for instance, \cite{NonlQCM}, are of the form of an ODE-IVP \eqref{DynSys}, i.e.,
\begin{align*}
	\begin{cases}
		\dot{x}(t) &= f(x(t),u(t),p), \\
		x(t_0) &= x_0, \\
	\end{cases}
\end{align*} 
with right-hand side 
\begin{align*}
	f(x,u,p) = 
	\left( \begin{array}{c}
		x_3 \\ x_4 \\ -\frac{p}{m_1} \left[ (x_1- x_2) + c (x_1 - x_2)^3 \right] - \frac{d_1}{m_1} (x_3 - x_4) \\ \frac{p}{m_2} \left[ (x_1 - x_2) + c (x_1 - x_2)^3 \right] + \frac{d_1}{m_2} (x_3 - x_4) -\frac{k_2}{m_2} (x_2 - u)
	\end{array} \right). 
\end{align*}
In addition, we choose the initial state $x_0 = (0,0,0,0)^\top$ for our example. The acceleration of the upper body mass $y$ is described by an output equation \eqref{OutputEqu}
of the form
\begin{align*}
	y(t) = h(x(t),p)
\end{align*}  
with right-hand side 
\begin{align*}
	h(x,p) = -\frac{p}{m_1} \left[ (x_1- x_2) + c (x_1 - x_2)^3 \right] - \frac{d_1}{m_1} (x_3 - x_4),
\end{align*}
independent of the road profile $u$.
For our example, the remaining model parameters (in SI-units) for the quarter-car are chosen to be
\begin{align}\label{ModelParameter}
	m_1 = 3600 ~[kg], ~ m_2 = 380 ~[kg], ~ k_2 = 1000 ~[kN/m], ~ d_1 = 34 ~[kN \cdot s/m], ~ c = 40.
\end{align}

\begin{figure}[h]
	\centering
	\includegraphics[width = \textwidth]{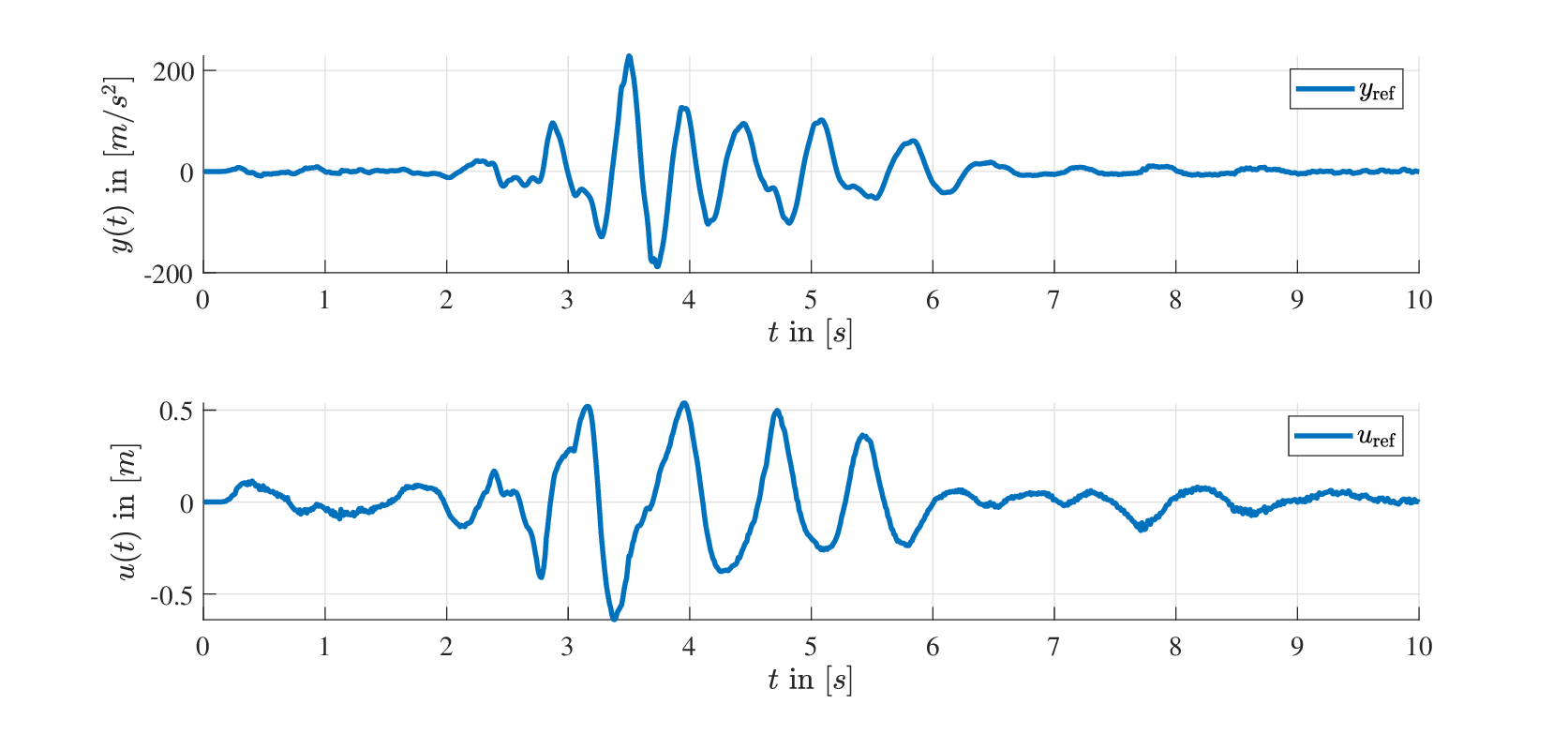} 
	\captionof{figure}{Reference output $y_\text{ref}$ (top) and input $u_\text{ref}$ (bottom) for a time interval of $I = [0, 10\text{s}]$.}\label{Figure Reference}
\end{figure}

In the forthcoming experiment, our focus is on the inverse problem, aiming to determine both the road profile $u$ and the stiffness parameter of the spring $p = k_1$ for a given quarter-car system. The objective is to manipulate these variables such that the acceleration of the upper body $y$ closely replicates a predefined reference output $y_\text{ref}$. 
For this simulation study, the provided reference acceleration $y_\text{ref}$ spans the time interval $I = [t_0, t_f] = [0, 10\text{s}]$ and is graphically depicted in Figure \ref{Figure Reference}.
\begin{remark}
	The reference output $y_\text{ref}$ is generated by a forward simulation of the quarter-car using a given road profile $u_\text{ref}$, also shown in Figure \ref{Figure Reference}, and the stiffness coefficient $p_\text{ref} = 230 ~[kN/m]$. Therefore, the results from solving the tracking problem with the suggested method can be compared to the reference solutions $u_\text{ref}$ and $p_\text{ref}$. Nevertheless, note that the objective function contains regularization terms such that one may not expect to exactly reproduce these input references.  
\end{remark}

\noindent
In formulating the objective function weights for Problem \ref{NonlinearODEOCP}, we chose $Q = 0.1$ and $T = 0.001$, with regularization parameters $\alpha_u = 30$ and $\alpha_p = 1e-10$. Additionally, we imposed box-constraints solely on the parameter, defined by the lower bound $p^\text{low} = 0.85 \cdot p^{(0)}$ and upper bound $p^\text{up} = 1.15 \cdot p^{(0)}$. From an application standpoint, the parameters computed during iterations, denoted as $p^{(k)}$, are permitted to deviate by at most $15\%$ from the initial value $p^{(0)}$.

To address the ODE optimal tracking control problem \ref{NonlinearODEOCP}, we have employed the Gauss-Newton method as presented in Algorithm \ref{Gauss-Newton for ODE Optimal Control}, incorporating the gradient descent method outlined in Algorithm \ref{GDMethod} for solving the auxiliary problems detailed in \ref{LQODEOCP}. Details regarding initial guesses and algorithmic parameters for the methods utilized can be found in Table \ref{Algorithmic Parameter}.

\begin{figure}[h]
	\centering
	\includegraphics[width = \textwidth]{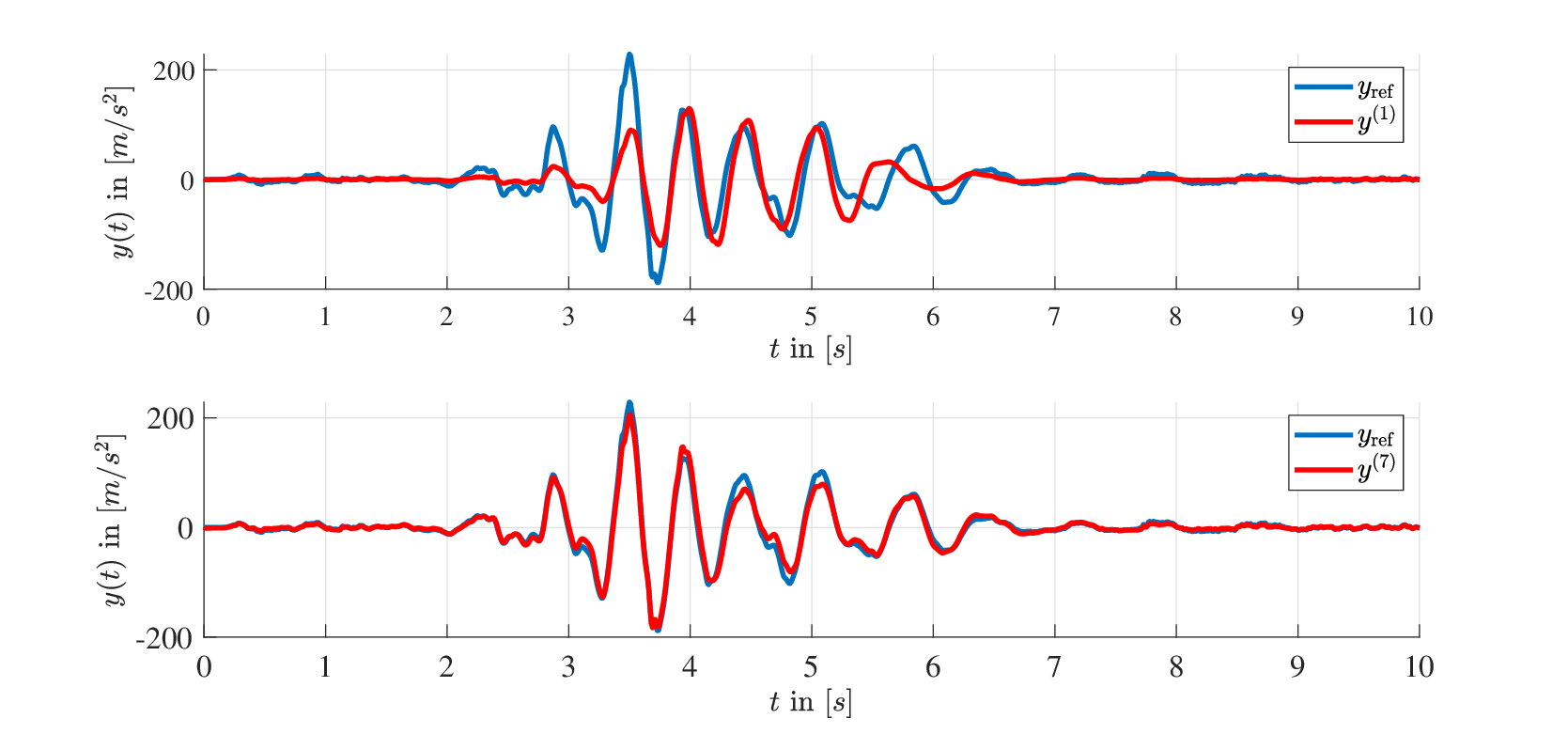}
	\captionof{figure}{Simulation $y^{(1)}$ (top) and $y^{(7)}$ (bottom), both compared to the reference $y_\text{ref}$.}\label{GDIteration17}
\end{figure}

\begin{table}[h]
	\centering
	\begin{tabular}{|m{0.25\textwidth}|m{0.25\textwidth}|}
		\hline
		Algorithmic parameter  & Algorithmic parameter \\
		for the Gauss-Newton method & for the gradient descent method \\
		\hline 
		\vspace{0.2cm} & \vspace{0.2cm}\\
		$u^{(0)}(t) = 0$ for $t \in [0, 10]$, & $\delta u^{(0)}(t) = 0$ for $t \in [0, 10]$, \\
		$p^{(0)} = 230$, & $\delta p^{(0)} = 20.7$,\\
		$\beta = 0.75$, & $\hat{\beta} = 0.3$, \\
		$\sigma = 1e-04$ & $\hat{\sigma} = 1e-04$, \\
		\hline
	\end{tabular} 
	\caption{Algorithmic parameters for the Gauss-Newton method and the gradient descent method.}\label{Algorithmic Parameter}
\end{table}
To solve the ODEs embedded within the algorithms, we decided for utilize the 'ode15s' solver from the MATLAB \cite{MATLAB} Toolbox. We established a time grid with a constant step size of $dt = 0.01\text{s}$. In Figure \ref{GDIteration17}, we display some of the output signals $y^{(k)}$ computed during the initial $N = 7$ iterations.

Notably, the final computed acceleration $y^{(7)}$ demonstrates a remarkable proximity to the reference signal $y_\text{ref}$, in stark contrast to the initial iterate $y^{(1)}$.\\

\noindent
The decreasing function values $J_k$ serve as quantifiable evidence of this outcome, with details provided in Table \ref{Function ValuesGD}. In our analysis, we have halted the iterations after $7$ runs, as additional iterations did not yield significant improvements in the results within the considered setup.\\

\begin{table}[h]
	\centering
	\begin{tabular}{|m{0.2\textwidth}|m{0.2\textwidth}|m{0.2\textwidth}|m{0.2\textwidth}|}
		\hline
		Iteration $k$  & Function value $J_k$ & Data misfit & Regularization terms\\
		\hline 
		$0$ & $1077.12$ & 1074.98 & 2.14\\ 
		$1$ & $403.68$ & 400.53 & 3.15 \\
		$3$ & $77.32$ & 73.12 & 4.20 \\
		$7$ & $35.20$ & 30.30 & 4.90 \\
		\hline
	\end{tabular}
	\caption{Function values $J_k$ for the iterations $k = 0,1,3,7$.}\label{Function ValuesGD}
\end{table} 

\begin{figure}[h]
	\centering
	\includegraphics[width = \textwidth]{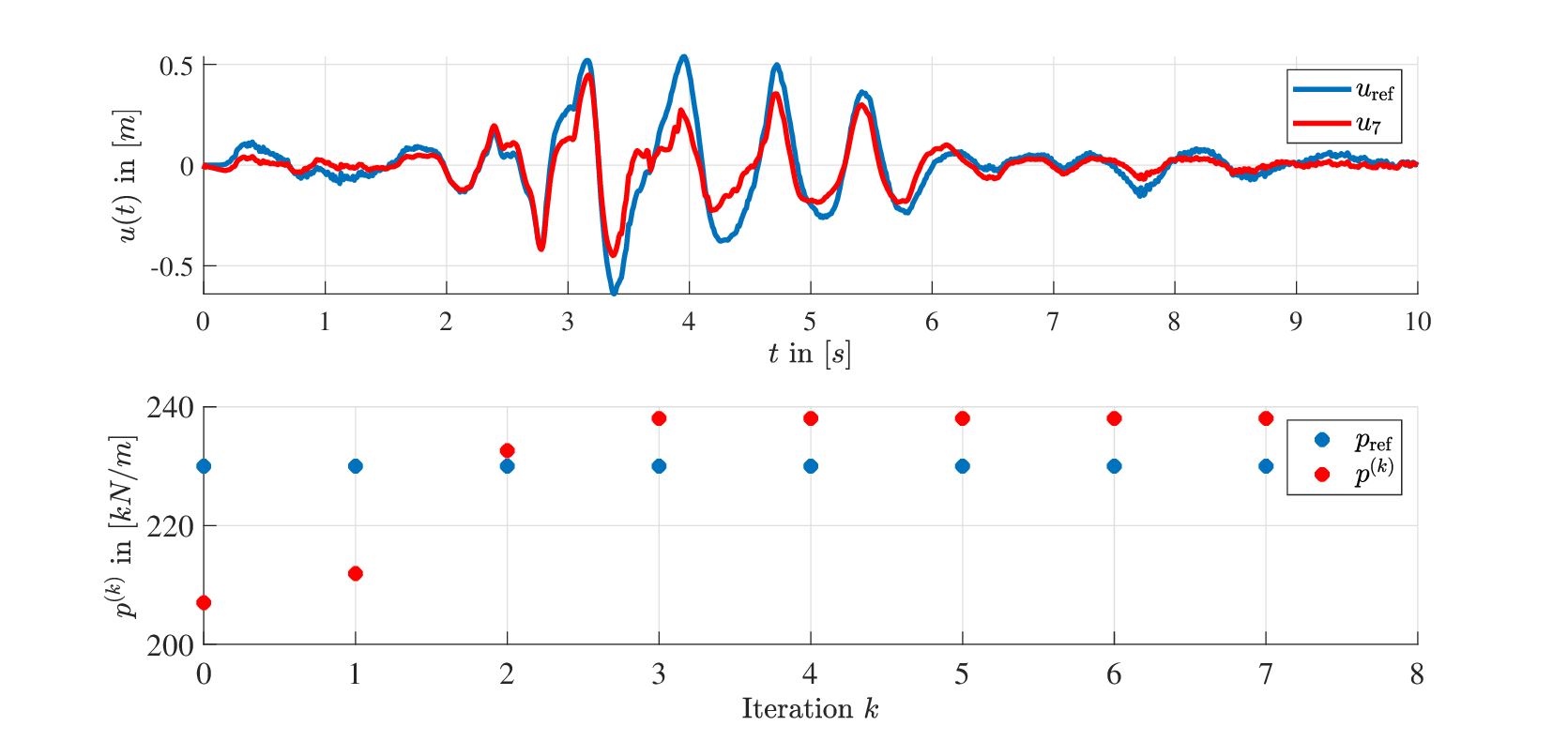} 
	\captionof{figure}{Computed road profile $u^{(7)}$ compared to $u_\text{ref}$ (top); computed parameters $p^{(k)}$, $k = 1,...,7$, compared to $p_\text{ref}$ (bottom).}\label{GDInputParameter}
\end{figure}

\noindent
Additionally, we have conducted a comparison between the computed road profile $u^{(7)}$ and parameter $p^{(7)}$ after the final iteration and their respective reference values $u_\text{ref}$ and $p_\text{ref}$, as depicted in Figure \ref{GDInputParameter}. Notably, the parameter values $p^{(k)}$ increased until reaching the upper bound of the box-constraints by the third iteration. While the computed road profile $u^{(7)}$ does not align well with the reference, the combination of this road profile with the parameter value led to a satisfactory outcome in the resulting output signal $y^{(7)}$.
\\

\begin{figure}[!h]
	\centering
	\includegraphics[width = \textwidth]{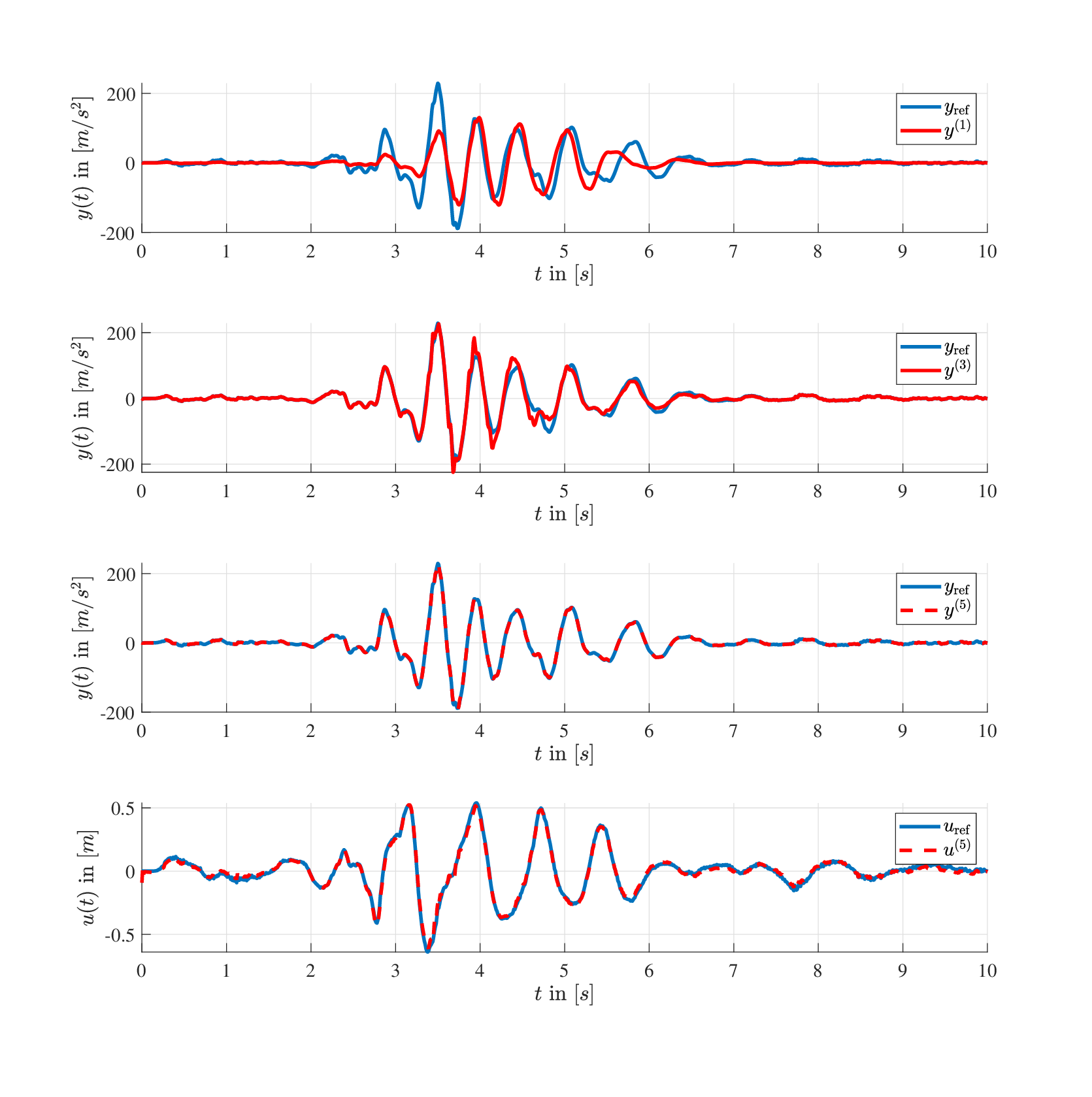} 
	\captionof{figure}{Simulations $y^{(1)}$ (first), $y^{(3)}$ (second) and $y^{(5)}$ (third) compared to reference output $y_\text{ref}$; computed road profile $u^{(5)}$ to the reference $u_\text{ref}$ (fourth).} \label{ResultsRiccati}
\end{figure}

\noindent
In a second experiment, we have assumed that the exact stiffness parameter $p$ corresponding to $y_\text{ref}$ is given, that is,
\begin{align*}
	p_\text{ref} = 230 ~[kN/m].
\end{align*}
Thus, we have fixed the parameter $\bar{p} = p_\text{ref}$ and have focused on the computation of the road profile $u$ for the reference acceleration $y_\text{ref}$ by applying the Gauss-Newton method for problem \eqref{ODEOptControl_u} including the Riccati-approach to solve the auxiliary problem \eqref{RiccatiAuxiliaryProblem_u} explicitly in each iteration step. The algorithmic parameters for the Gauss-Newton method have been chosen as in the previous simulation study, i.e.,
\begin{align}
	u^{(0)}(t) = 0,~ t \in [0, 10], \quad \beta = 0.75, \quad \sigma = 1e-04.
\end{align}
In total, we have run $N = 5$ iterations and have set the regularization parameter $\alpha_u = 30$. 
\begin{remark}
	Since $h_u(x,\bar{p}) = 0$ for the considered quarter-car and $n_u = 1$, we have that $R(t) = \alpha_u$ for all $t \in [0, 10]$. Therefore, the regularization parameter $\alpha_u$ had to be chosen larger than $0$ to guarantee that the matrix $R(t)$ is positive definite. Furthermore, for the matrix $\hat{Q}(t)$ in \eqref{HatQ} we have that
	\begin{align*}
		\hat{Q}(t) = \left( \begin{array}{cc}
			C(t)^\top Q C(t) & -C(t)^\top Q r(t) \\
			-r(t)^\top Q C(t) & r(t)^\top Q r(t) - k(t)^\top R(t)^{-1} k(t)
		\end{array}\right),
	\end{align*}
	with $Q = 0.1 > 0$, and, therefore, $\hat{Q}(t)$ is positive semi-definite for sufficiently large chosen $\alpha_u$. In addition, since also $T = 0.001 > 0$, the matrix $\hat{T}$ in \eqref{HatT} is positive semi-definite. Hence, all assumptions for the Riccati-approach are satisfied.
\end{remark}
\noindent
The numerical results for this example are presented in Figure \ref{ResultsRiccati}. Following the first iteration, an evident phase shift of the simulation $y^{(1)}$ within the interval $[4.5\text{s}, 6.5\text{s}]$ in comparison to the reference $y_\text{ref}$ is noticeable. Furthermore, there is a discrepancy in the amplitudes of both signals. Subsequent to the third iteration, the phase shift in the simulated output signal $y^{(3)}$ was rectified, and after two additional iterations, the amplitudes of $y^{(5)}$ aligned closely with those of the reference signal $y_\text{ref}$.

\noindent
Figure \ref{ResultsRiccati} also shows the computed road profile $u^{(5)}$ in comparison to the reference profile $u_\text{ref}$, which obviously seem to be very close. We further want to quantify the difference of the computed output signals to the given reference by considering the function values $J_k$ during the iteration process and, particularly, the data misfit given by 
\begin{equation}\label{DataMisfit}
	\frac{1}{2} \left\Vert Q^{\frac{1}{2}} \left[\S(u^{(k)},p_\text{ref}) - y_\text{ref} \right] \right\Vert_{L^\infty(I;\R^{n_y})}^2  + \frac{1}{2} \left\Vert T^{\frac{1}{2}} \left[\S(u^{(k)},p_\text{ref})(t_f) - y_\text{ref}(t_f) \right] \right\Vert_{\R^{n_y}}^2.
\end{equation} 
The corresponding values are stated in Table \ref{Function ValuesGN}.
\begin{table}[h]
	\centering
	\begin{tabular}{|m{0.2\textwidth}|m{0.2\textwidth}|m{0.2\textwidth}|m{0.2\textwidth}|}
		\hline
		Iteration $k$  & Function value $J_k$ & Data misfit \eqref{DataMisfit} & Regularization term\\		   
		\hline 
		$0$ & $1077.28$ & $1077.28$ &  $0$ \\
		$1$ & $415.16$ & $414.32$ & $0.84$\\
		$3$ & $57.33$ & $53.29$ & $4.04$\\
		$5$ & $4.73$ & $0.49$ & $4.24$\\
		\hline
	\end{tabular}
	\caption{Function values $J_k$ for the iterations $k = 0,1,3,5$.}\label{Function ValuesGN}
\end{table}
A rapid decrease in the function values $J_k$ and the data misfit is evident within the initial five iterations. Consequently, in this example, our presented method demonstrates the capability to achieve excellent results with only a few iterations.

\begin{figure}[ht]
	\centering
	\includegraphics[width = \textwidth]{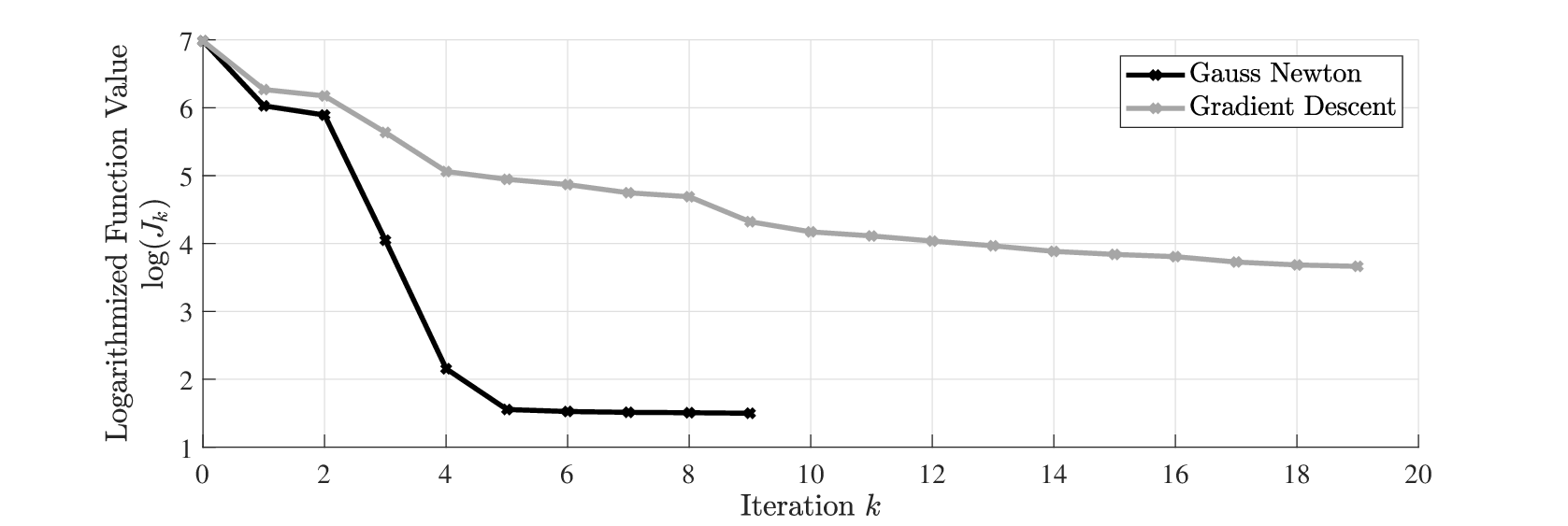}
	\captionof{figure}{Comparison of the (logarithmic) function values in the Gauss-Newton method and the gradient descent method.}\label{ComparisonGNvsGD}
\end{figure}

Furthermore, we have conducted a comparison between the outcomes obtained through the proposed Gauss-Newton approach and those achieved by directly employing the function space gradient descent method to address problem \ref{ODEOptControl_u}. Figure \ref{ComparisonGNvsGD} illustrates the resulting function values $J_k$ for both methods across the iteration process. It is evident that the proposed Gauss-Newton method converges much more rapidly than the gradient descent function space method and requires significantly fewer iterations to attain very satisfactory results in this example. 

\section{Summary \& Conclusions}\label{sec5}
We introduced and analyzed a method tailored to address ODE optimal control problems featuring tracking objectives and box-constraints within a specific function space framework. By treating the problem as an infinite-dimensional optimization challenge, we apply the Gauss-Newton method due to its suitability for nonlinear least squares problems, leveraging a projection technique to incorporate the box-constraints. Each iteration of our method involves solving a linear-quadratic ODE optimal control problem. We propose two distinct approaches from ODE optimal control to address these subproblems. Additionally, we provide a global convergence result for the Gauss-Newton method based on optimization theory in $L^2$. Finally, we present a specific application scenario where we apply the method using both approaches for the subproblems and discuss the results.

Compared to discretization approaches, the function space method offers several advantages. It avoids the curse of dimensionality and relies solely on integration techniques for solving ODEs. Moreover, discretization occurs at a later stage, enhancing flexibility. Sensitivity to algorithmic parameters is lower due to the reduced parameter set, resulting in increased robustness. On the other hand, complex constraints beyond box-constraints pose challenges in function space methods and often necessitate alternative techniques. Discretization approaches are more suitable in such cases. Practical applications often involve other forms of constraints, diminishing the utility of function space methods. If control and model parameters are optimized simultaneously, the advantages of function space methods diminish as the auxiliary problem can no longer be solved explicitly. Choosing algorithmic parameters becomes more challenging in such scenarios, making discretization approaches a viable alternative.

Especially in applications in the field of numerical system simulation, the model equations for the dynamical system are often build within a commercial software tool and one may have limited access, for instance, due tue their complexity. For the use of the presented method one not necessarily requires the full model equations. It is sufficient to have evaluations of the right-hand sides as well as their Jacobians provided by the tool.

\section*{Acknowledgments}
Claudia Schillings acknowledges support from MATH+ project EF1-19: Machine Learning Enhanced Filtering Methods for Inverse Problems, funded by the Deutsche Forschungsgemeinschaft (DFG, German Research
Foundation) under Germany's Excellence Strategy - The Berlin Mathematics
Research Center MATH+ (EXC-2046/1, project ID: 390685689). 

\section*{Appendix}
\textit{Proof of Theorem \ref{Theorem_Adjoint}}.\label{Proof_Adjoint} We show that, for given $\hat{u}$ with $\tilde{\S}(\hat{u}) = \hat{x}$, the following equation holds:
\begin{align*}
	\left\langle S^\prime (\hat{u})(\delta u),\delta y  \right\rangle_{L^2(I;\R^{n_y})} = \left\langle (\delta u), \left( S^\prime (\hat{u})\right)^* (\delta y) \right \rangle_{L^2(I;\R^{n_u})}.
\end{align*}
Let $\phi(t,t_0)$ denote the state transition matrix that is the unique global solution of the matrix differential equation initial value problem
\[ \frac{\partial}{\partial t} \phi(t,t_0) = f_x[t]^\top \phi(t,t_0), \quad  \phi(t_0,t_0) = I_{n_x \times n_x}.\]
Then, the solution of system corresponding to $S^\prime (\hat{u})$ is given by 
\begin{align*}
	\delta y(t) = h_x[t] \int_{t_0}^t \phi(t,s) f_u[s]\delta u(s)   ~ds + h_u[t] \delta u(t) 
\end{align*}
for $t \in [t_0,t_f]$. We further compute
\begin{align*}
	\left\langle S^\prime (\hat{u})(\delta u), \delta y  \right\rangle_{L^2(I;\R^{n_y})}  &= \int_{t_0}^{t_f} \left( h_x[t] \int_{t_0}^t \phi(t,s)f_u[s]\delta u(s) ~ds + h_u[t] \delta u(t)\right)^\top \delta y(t) ~dt \\
	&=	\int_{t_0}^{t_f} \int_{t_0}^t \left[ h_x[t] \phi(t,s) f_u[s]\delta u(s) \right]^\top \delta y(t) ~ds~dt +\int_{t_0}^{t_f} \delta u(t)^\top h_u[t]^\top \delta y(t) ~dt.
\end{align*}
First, we focus on the first term for which we have
\begin{align*}
	\int_{t_0}^{t_f} \int_{t_0}^t \left[ h_x[t] \phi(t,s) f_u[s]\delta u(s) \right]^\top \delta y(t) ~ds~dt 	&= \int_{t_0}^{t_f} \int_{t_0}^t \delta u(s)^\top \left[ f_u[s]^\top \phi(t,s)^\top h_x[t]^\top \delta y(t) \right] ~ds~dt \\
	&= \int_{t_0}^{t_f} \int_{s}^{t_f} \delta u(s)^\top \left[ f_u[s]^\top \phi(t,s)^\top h_x[t]^\top \delta y(t) \right] ~dt~ds \\
	&= \int_{t_0}^{t_f} \delta u(s)^\top \left[f_u[s]^\top \int_{t_f}^{s} \phi(t,s)^\top (-h_x[t]^\top) \delta y(t) ~dt \right] ~ds \\
	&= \left \langle \delta u, f_u[\cdot]^\top \int_{t_f}^{\cdot} \phi(t,\cdot)^\top (-h_x[t]^\top)  \delta y(t) ~dt \right \rangle_{L^2(I;\R^{n_u})}.
\end{align*}
Note, that we applied Fubini in the second equality. For the second term we have
\begin{align*}
	\int_{t_0}^{t_f} \delta u(t)^\top h_u[t]^\top \delta y(t) ~dt = \left \langle \delta u, h_u[\cdot]^\top \delta y \right \rangle_{L^2(I;\R^{n_u})}.
\end{align*}
Finally, we get
\begin{align*}
	\left\langle S^\prime (\hat{u})(\delta u), \delta y \right\rangle_{L^2(I;\R^{n_y})}  &= \left \langle \delta u, f_u[\cdot]^\top \int_{t_f}^{\cdot} \phi(t,\cdot)^\top (-h_x[t]^\top) \delta y(t) ~dt \right \rangle_{L^2(I;\R^{n_u})} + \left \langle \delta u, h_u[\cdot]^\top \delta y\right \rangle_{L^2(I;\R^{n_u})} \\
	&= \left \langle \delta u, f_u[\cdot]^\top \int_{t_f}^{\cdot} \phi(t,\cdot)^\top (-h_x[t]^\top) \delta y(t) ~dt + h_u[\cdot]^\top \delta y \right \rangle_{L^2(I;\R^{n_u})} \\
	&= \left\langle \delta u, \left( S^\prime (\hat{u})\right)^\ast (\delta y) \right \rangle_{L^2(I;\R^{n_u}) },
\end{align*}
whereas $\left(S^\prime (\hat{u})\right)^\ast$ is the solution operator of the given system from the theorem and we used the fact that the solution of the ODE end-value problem \eqref{AdjDynSys1} is given by
\begin{align*}
	\delta x(t) = \int_{t_f}^{t} \phi(s,t)^\top(-h_x[s]^\top)\delta y(s) ~ds.
\end{align*}
\qed

\bibliographystyle{alpha}
\bibliography{A_Gauss_Newton_Method_for_ODE_Optimal_Tracking_Control}

\clearpage

\end{document}